\documentclass[journal,twoside,web]{ieeecolor}
\usepackage{generic}
\usepackage{cite}
\usepackage{amsmath,amssymb,amsfonts}
\usepackage{algorithmic}
\usepackage{graphicx}
\usepackage{algorithm,algorithmic}
\usepackage{hyperref}
\usepackage{bm}
\usepackage{mathrsfs}
\usepackage{mathtools}
\usepackage{textcomp}
\usepackage{xcolor}
\usepackage{lineno,hyperref,float,amsfonts,epsfig,graphics,booktabs,multirow,stfloats}
\usepackage{cite}

\usepackage{amsthm}
\usepackage{comment}

\newcommand{\tabincell}[2]{\begin{tabular}{@{}#1@{}}#2\end{tabular}}

\theoremstyle{definition}

\theoremstyle{plain}
\newtheorem{myTheo}{Theorem}

\newtheorem{myLemma}{Lemma}
\theoremstyle{remark}
\newtheorem{myRemark}{Remark}

\newtheorem{myAssump}{Assumption}

\DeclareMathOperator{\tr}{tr}
\DeclareMathOperator{\diag}{diag}
\DeclareMathOperator{\rank}{rank}

\hypersetup{hidelinks=true}
\usepackage{textcomp}
\def\BibTeX{{\rm B\kern-.05em{\sc i\kern-.025em b}\kern-.08em
    T\kern-.1667em\lower.7ex\hbox{E}\kern-.125emX}}
\begin{document}
\title{Adaptive Robust Control for Uncertain Systems with Ellipsoid-Set Learning}
\author{Xuehui Ma, Shiliang Zhang,~\IEEEmembership{Member,~IEEE}, Zhiyong Sun,~\IEEEmembership{Member,~IEEE}, \\Xiaohui Zhang, Sabita Maharjan,~\IEEEmembership{Senior Member,~IEEE}
\thanks{This work was partially supported by the National Natural Science Foundation of China under Grant 62073259. \\Xuehui Ma and Xiaohui Zhang are with the School of Automation and Information Engineering, Xi'an University of Technology, Xi'an, China (xuehui.yx@gmail.com, xiaohui.xaut@gmail.com). Shiliang Zhang and Sabita Maharjan are with the Department of Informatics, University of Oslo, Oslo, Norway (\{shilianz, sabita\}@ifi.uio.no). Zhiyong Sun is with College of Engineering, Peking University, Beijing, China (zhiyong.sun@pku.edu.cn). \\© 2026 IEEE.  Personal use of this material is permitted.  Permission from IEEE must be obtained for all other uses, in any current or future media, including reprinting/republishing this material for advertising or promotional purposes, creating new collective works, for resale or redistribution to servers or lists, or reuse of any copyrighted component of this work in other works.}
}

\maketitle

\begin{abstract}

Despite the celebrated success of stochastic control approaches for uncertain systems, such approaches are limited in the ability to handle non-Gaussian uncertainties. This work presents an adaptive robust control for linear uncertain systems, whose process noise, observation noise, and system states are depicted by ellipsoid sets rather than Gaussian distributions. We design an ellipsoid-set learning method to estimate the boundaries of state sets, and incorporate the learned sets into the control law derivation to reduce conservativeness in robust control. Further, we consider the parametric uncertainties in state-space matrices. Particularly, we assign finite candidates for the uncertain parameters, and construct a bank of candidate-conditional robust control problems for each candidate. We derive the final control law by aggregating the candidate-conditional control laws. In this way, we separate the control scheme into parallel robust controls, decoupling the learning and control, which otherwise renders the control unattainable. We demonstrate the effectiveness of the proposed control in numerical simulations in the cases of linear quadratic regulation and tracking control.

\end{abstract}


\section{Introduction}\label{section1}
This work looks into the control of uncertain linear systems interrupted by unknown noises, an extensively investigated topic yet still with challenging issues in mitigating the negative impact of uncertainties. Uncertainties inevitably exist in the control of practical systems due to unknown or varying system dynamics, disturbances, and noises, which complicate the maintenance of control over time. To address this challenge, stochastic adaptive control~\cite{aastrom2012introduction,åström2008adaptive} has emerged and served as an established method in handling uncertainties, which utilizes noisy data to learn the system uncertainties and eliminate their impact on the control. Notably, Åström and Wittenmark~\cite{ASTROM1973185} proposed the self-tuning regulators for unknown systems described by controlled autoregressive moving average models (ARMA). They estimate the unknown parameters using recursive least squares (RLS) and derive the control law via minimum variance control (MVC). Deshpande \textit{et al.}~\cite{DESHPANDE1973107} represented stochastic systems by state-space equations, and utilized the \textit{a-posteriori} probabilities to learn the unknown parameters. They integrated the learned parameters with linear quadratic Gaussian control (LQG) to derive an adaptive control law for stochastic systems. With their recognized advantages in learning and reducing uncertainties, the mentioned classical methods have been extended to adaptive model predictive control~\cite{mesbah2018stochastic}, adaptive dual control~\cite{filatov2004adaptive,tse1973wide}, anti-disturbance control~\cite{narendra2001stochastic,ma2022active}, and fault-tolerant control~\cite{liu2022dual,liu2021reliable}. However, these works assume a prior probability distribution for the studied uncertainties, usually a Gaussian distribution. This assumption is often hampered by a shortage of information to build the required probabilistic models~\cite{becis2008state}, \textit{e.g.}, for the probability distribution of the unknown noises and parameters.

Different from stochastic adaptive control that requires explicit probability distributions to describe the uncertainties, robust control quantifies the uncertainties by bounded sets~\cite{zhou1998essentials,lee1997worst}. The rationale behind is that, in most practical scenarios, the boundaries of the possible values for the uncertainties, \textit{e.g.}, system noises, model parameters, or system states, are relatively easy to obtain by examining the physics of the controlled system~\cite{yang2009set,maksarov1996state}. Such boundaries formulate the bounded sets in robust control, which then aims to develop a conservative control policy that performs the best for the worst uncertainties within the bounded sets. The representative example of robust control is the min-max control method. \textit{E.g.}, Kothare \textit{et al.}~\cite{kothare1996robust} presented a min-max model predictive control (MPC) for uncertain systems with unknown parameters confined in polytopic sets. Bertsimas and Brown~\cite{bertsimas2007constrained} used the robust optimization framework to solve the min-max control law for the linear quadratic control problem, where the system states and noises belong to ellipsoid sets. Nevertheless, robust control assumes that the uncertainty sets are fixed over the control horizon. Such fixed sets cannot reflect the actual system dynamics over time, resulting in an inadequate representation of the uncertainties during the control. Furthermore, these fixed sets have to be large enough to cover all potential worst-case scenarios, leading to an over-conservative control policy that compromises control accuracy.

Addressing the over-conservativeness issue in robust control necessitates narrowing the range of uncertainties. Studies have been conducted to reduce the size of bounded sets for uncertainties, via learning from system observations to gain knowledge about the system dynamics and update the bounded sets. Set-membership filter has served as a powerful tool to learn and reduce the bounded set representing the uncertainties~\cite{LI2023110842,LI2023110874,LIU2021109684}, and has drawn attention in studies that incorporate bounded set learning in robust control. Relevant work includes Lorenzen \textit{et al.}~\cite{lorenzen2019robust}, Köhler \textit{et al.}~\cite{kohler2021robust}, Lu \textit{et al.}~\cite{8814456}, Sasfi \textit{et al.}~\cite{sasfi2023robust}, Tranos \textit{et al.}~\cite{tranos2023self}, Parsi \textit{et al.}~\cite{9780027,parsi2020active,parsi2023dual}, and Arcari \textit{et al.}~\cite{10179989} that applied the set-membership filter to update the bounded sets of system parameters and states in robust model predictive control. Their simulation results demonstrate significant enhancement in control performance. However, a notable drawback of set-membership filters is their reliance on polytope sets for uncertainty learning~\cite{lorenzen2019robust,maksarov1996state}. Polytope sets can lead to a dramatic increase in algorithm complexity as the system dimension grows, which jeopardizes the utility of set-membership filters in practical applications. 

An emerging approach to reducing computational complexity in bounded set learning is to use ellipsoid sets for uncertainty representation~\cite{maksarov1996state,TANG202013688}. Several studies~\cite{iannelli2020multiobjective,iannelli2020structured,qian2010adaptive,parsi2022scalable,sasfi2023robust,ma2024adaptive} have utilized ellipsoid sets to represent uncertain parameters in robust control. In these studies, the ellipsoid-set learning method is applied to update the parameter sets, which are then incorporated into the design of adaptive robust control strategies. Nevertheless, these methods assume the system states are directly measurable. This is a strong assumption that often does not hold in real-world applications. This limitation significantly restricts the applicability of the existing ellipsoid-set-based adaptive robust control methods, particularly for systems whose states are not directly measurable. Recently, studies by Velázquez \textit{et al.}~\cite{velazquez2021attractive}, Wang  \textit{et al.}~\cite{wang2022ellipsoidal}, and Ma \textit{et al.}~\cite{ma2024robust} introduced ellipsoid set-membership filter to estimate the state in deriving robust control law. However, their methods only apply to control systems whose parameters are determined and known. That is, their developed methods require that the parameters of the system model are completely known, an assumption that hardly holds in practice. 

In this work, we emphasize both the state estimation and parameter learning when designing robust control for uncertain systems. Furthermore, we aim to develop a strategy that is free from uncertainty distribution requirements or strict assumptions of directly measurable states. We present an ellipsoid-set learning algorithm to estimate the ellipsoid-set of system states, and utilize the updated set to derive the min-max control law based on robust optimization. As for the uncertain parameters in state-space matrices, we note that integrating the uncertain parameter learning into the min-max robust control leads to the coupling between learning and control. Such coupling can render the solution to adaptive control law derivation unattainable~\cite{iannelli2020multiobjective,ma2024adaptive}. This issue of coupled learning and control is also known as the dual control problem~\cite{wittenmark1995adaptive,filatov2000survey,MA2022157}. Though a number of sub-optimal solutions have been developed for dual control, these available solutions are not applicable to solve issues in the ellipsoid-set-based adaptive robust control~\cite{hewing2020learning}.

In this work, we address the issue of learning-control coupling in ellipsoid-set-based adaptive robust control. We design a finite candidate set to approximate the unknown parameters, and utilize \textit{a-posteriori} probability to learn the weight for each of the individual candidate values. For each candidate value of uncertain parameters, we formulate a candidate-conditional robust control problem and derive the control law. In this way, we obtain separated and decoupled parameter learning and robust control that can otherwise compromise the uncertain system control. We derive our final adaptive robust control law by aggregating the candidate-conditional robust control laws weighted by their \textit{a-posteriori} probabilities. 
We conduct numerical simulations to demonstrate our approach in comparison with min-max control with fixed sets~\cite{bertsimas2007constrained} and optimal robust control with known system parameters~\cite{ma2024robust}. We summarize the main contributions of this study as follows.
\begin{enumerate}
    \item We propose to use ellipsoid sets to represent system states, process noises, and observation noises in our state space model. This deviates from conventional Markov models, which typically assume Gaussian distributions for uncertainties. The adopted ellipsoid set provides a more flexible and robust representation of the uncertainties, especially for systems with bounded uncertainties whose statistics and distributions are difficult to obtain. 
    \item We develop an ellipsoid-set learning algorithm to estimate the system states. This differs from existing approaches relying on Kalman filter, which are inherently limited to Gaussian noise assumptions. The developed algorithm is able to handle bounded noises with unknown distributions and enhance system state estimation in \textit{non-Gaussian} conditions.
    \item Most adaptive control strategies directly learn unknown parameters for uncertain systems. In our study, we propose candidate values for unknown parameters and derive the learning of weights for the candidate values. This naturally decouples the parameter learning from control design, which otherwise significantly complicates the control law derivation.
\end{enumerate}

The rest of the paper is organized as follows. Section \ref{section2} formulates the linear quadratic control problem for uncertain systems with process and observation noises bounded by ellipsoid sets. Section \ref{section3} presents the derivation of the adaptive robust control law for linear quadratic regulation and tracking control problems. In Section \ref{section5}, we conduct numerical simulations to demonstrate the effectiveness of our methods. We conclude this work in Section~\ref{section6}.\\

\noindent\textit{Notation:} \\
In this work, we define the ellipsoid with centre $\bm{c}$ and shape matrix $\bm{P}$ by $E(\bm{c},\bm{P})\equiv \{\bm{x}\in \mathbb{R}^n|(\bm{x}-\bm{c})^T\bm{P}^{-1}(\bm{x}-\bm{c})\leq 1\}$, where matrix $\bm{P}$ is symmetric positive definite. $\bm{I}$ and $\bm{0}$ denote the identity matrix and zero matrix, respectively, with appropriate dimensions. 
The Minkowski sum of two ellipsoid sets $E(\bm{c}_1,\bm{P}_1)$ and $E(\bm{c}_2,\bm{P}_2)$ is defined as adding each vector in $E(\bm{c}_1,\bm{P}_1)$ to each vector in $E(\bm{c}_2,\bm{P}_2)$: $E(\bm{c}_1,\bm{P}_1) \oplus E(\bm{c}_2,\bm{P}_2)=\{ \bm{x}_1+\bm{x}_2| \bm{x}_1 \in E(\bm{c}_1,\bm{P}_1), \bm{x}_2 \in E(\bm{c}_2,\bm{P}_2) \} $.

\section{Problem Statement} \label{section2}
We consider the control of a discrete-time uncertain system with ellipsoid-set bounded noises. We describe the system by a state-space model shown as 
\begin{equation}
	\begin{aligned}
		\bm{x}(k+1)=\bm{A}(k,\bm{\theta})\bm{x}(k)+\bm{B}(k,\bm{\theta})\bm{u}(k)+\bm{w}(k),
	\end{aligned}
\end{equation}
\begin{equation}
	\begin{aligned}
		\bm{y}(k) = \bm{C}(k,\bm{\theta})\bm{x}(k)+\bm{v}(k),
	\end{aligned}
\end{equation}
where $\bm{x}(k) \in \mathbb{R}^n$ is the state vector, $\bm{u}(k) \in \mathbb{R}^r$ is the control vector, $\bm{y}(k) \in \mathbb{R}^m$ is the system observation. $\bm{w}(k) \in \mathbb{R}^n$ and $\bm{v}(k) \in \mathbb{R}^m$ are the process noise and observation noise, respectively.
We make the following assumptions on the system states and noises.
\begin{myAssump}\label{Assump_obser}
The states of the system (1) are not directly measurable, but are reconstructible through the system dynamics and output measurements.
\end{myAssump}
\begin{myAssump}\label{Assump_noise}
The process and observation noises are assumed to be bounded by the following ellipsoid sets
\begin{equation}\label{noise_w}
	\begin{aligned}
		\bm{w}(k) \in E(\bm{0},\bm{P}_{w}(k)),
	\end{aligned}
\end{equation}
\begin{equation}\label{noise_v}
	\begin{aligned}
		\bm{v}(k) \in E(\bm{0},\bm{P}_{v}(k)) ,
	\end{aligned}
\end{equation}
where $\bm{P}_{w}(k)\in \mathbb{R}^{n \times n} $ and $\bm{P}_{v}(k)\in \mathbb{R}^{m \times m}$ are the matrices defining the shape, size, and orientation of the ellipsoids.
\end{myAssump}
\begin{myAssump} \label{Assump_state0}
The initial state vector $\bm{x}(0)$ is unknown and bounded by the following ellipsoid set:
\begin{equation}\label{state_x0}
	\begin{aligned}
		\bm{x}(0)\in E(\bar{\bm{x}}_0,\bm{P}_0) ,
	\end{aligned}
\end{equation}
where $\bar{\bm{x}}_0$ is the center of the ellipsoidal set, and $\bm{P}_0$ is the ellipsoid shape matrix. The initial center and shape matrix can be selected based on system characteristics and prior knowledge of the system of interest.    
\end{myAssump}

Note that the value of state matrix $\bm{A}(k,\bm{\theta}) \in \mathbb{R}^{n\times n}$, control matrix $\bm{B}(k,\bm{\theta})\in \mathbb{R}^{n\times r}$, and observation matrix $\bm{C}(k,\bm{\theta})\in \mathbb{R}^{m\times n}$ are all dependent on an unknown parameter vector $\bm{\theta}\in \mathbb{R}^{l}$. We define a parameter space $\bm{\Theta}$ - namely the range of the uncertain parameter vector $\bm{\theta}$ - from which $\bm{\theta}$ takes its value, shown as
\begin{equation}\label{Theta_set}
	\begin{aligned}
		\bm{\Theta} = \{\bm{\theta}_1, \bm{\theta}_2, \cdots, \bm{\theta}_s\}.
	\end{aligned}
\end{equation}
The discrete parameter space appears only in a few practical scenarios, while the range of unknown parameters is continuous in most practical applications. We can discretize the continuous space as \eqref{Theta_set} to facilitate the decoupling between the parameter learning and the control, which leads to tractable control law derivation. Note that discretization of continuous space has been widely used in the development of estimation methods and adaptive control studies \cite{DESHPANDE1973107,shahab2021asymptotic,katsikas2001genetically}. We explain the parameter discretization in \textit{Remark} \ref{remark4}. 

%Such that we have an explicit range from which we can recognize the value of $\bm{\theta}$, regardless of whether the parameter space $\bm{\Theta}$ is discrete or continuous. 

In this work, we aim to gain the control sequence $\{{u}(k)\} \vert_{k=0}^{N-1}$, which minimizes the quadratic cost function 
\begin{equation} \label{J_index}
	\begin{aligned}
		J =  &\bm{x}^T(N)\bm{Q}(N)\bm{x}(N) \\
        &+\sum_{k= 0}^{N-1} \left[ \bm{x}^T(k)\bm{Q}(k)\bm{x}(k) +\bm{u}^T(k)\bm{R}(k)\bm{u}(k) \right],
	\end{aligned}
\end{equation}
where $\bm{Q}(k)\in\mathbb{R}^{n\times n} $ and $\bm{R}(k)\in\mathbb{R}^{r\times r}$ are positive semi-definite and positive definite symmetric matrices, respectively.
We formulate the following min-max robust control problem with the ellipsoid sets of process noise, observation noise, and initial state as constraints:
\begin{equation}\label{prob1}
	\begin{aligned}
		(\mathcal{P}):  &\min_{\{{u}(k)\} \vert_{k=0}^{N-1}} \max_{ \{\bm{w}(k), \bm{v}(k), \bm{x}(0) \} } J \\
		 s.t.  \quad & \bm{x}(k+1)=\bm{A}(k,\bm{\theta})\bm{x}(k)+\bm{B}(k,\bm{\theta})\bm{u}(k)+\bm{w}(k),\\
		\quad \quad \quad \quad & \bm{y}(k) = \bm{C}(k,\bm{\theta})\bm{x}(k)+\bm{v}(k),  k = 0,1, \cdots, N-1 .\\
		\quad \quad \quad \quad & \bm{w}(k) \in E(\bm{0},\bm{P}_w(k)), \\
		\quad \quad \quad \quad & \bm{v}(k) \in E(\bm{0},\bm{P}_v(k)), \\
		\quad \quad \quad \quad & \bm{x}(0)\in E(\bar{\bm{x}}_0,\bm{P}_0),\\
		\quad \quad \quad \quad & \bm{\theta} \in \bm{\Theta}= \{\bm{\theta}_1, \bm{\theta}_2, \cdots, \bm{\theta}_s\} , \\
	\end{aligned}
\end{equation}
The notes and challenges of our solutions to this control problem are described in the remarks as follows: 

\begin{myRemark}\label{remark4}
This remark provides a discretization approach for the continuous parameter space. We consider a one-dimensional parameter space $\bm{\Theta} = [\gamma_l, \gamma_u]$, where the parameters $\gamma_l$ and $\gamma_u$ denote the lower and upper boundaries, respectively. We divide this continuous space into $s$ sub-intervals uniformly, shown as $\bm{\Theta} = [\gamma_l, \gamma_u]= \bigcup_{i=1}^s[\gamma_i, \gamma_{i+1}]$. 
We choose the center $\bm{\theta}_i = (\gamma_i+\gamma_{i+1})/2$ of sub-interval to approximate the true parameter $\bm{\theta}$ that belongs to the sub-interval $\bm{\theta}\in[\gamma_i, \gamma_{i+1}]$. 
Then, the continuous space can be discretized as the discrete parameter space \eqref{Theta_set}. The discretization number is chosen by $s= \lceil \frac{\gamma_u-\gamma_l}{2\varepsilon} \rceil$, where $\varepsilon$ is the admissible error between the true parameter $\bm{\theta}$ and the approximated parameter $\bm{\theta}_i$, $\lceil x \rceil$ represents the least integer that is equal to or greater than $x$.
\end{myRemark}

\begin{myRemark}\label{remark1}
In problem $(\mathcal{P})$, we suppose that the system states are not directly measurable, and system noises are assumed to be ellipsoid-set-bounded rather than with explicit distributions, \textit{e.g.}, white Gaussian noises. Thus, the established state estimator like Kalman filter, where Gaussian noise is assumed, no longer fits for estimating the states in this case. Therefore, we are motivated to look into an iterative state bounding algorithm that can handle state unmeasurability excited by unknown noises bounded in ellipsoid sets, which can lead to the derivation of min-max robust controllers. 
\end{myRemark}

\begin{myRemark}\label{remark2}
It is possible to find a sequence of control signals $\{{u}(k)\} \vert_{k=0}^{N-1}$ that performs the best for the worst initial states and noises within the uncertainty sets described by ellipsoids. Nevertheless, the worst case might hardly happen in practical scenarios, which may render the obtained control policy to be over-conservative, especially when the ellipsoid set for the initial states is very large. This motivates us to design an approach that can reduce the conservativeness and improve the control performance. 
\end{myRemark}

\begin{myRemark}\label{remark3}
Deriving the adaptive robust control law requires the identification of unknown parameter $\bm{\theta}$, while the control law derivation affects the identified result of $\bm{\theta}$. Such a coupling between the control and the unknown parameter learning is an issue commonly observed in dual control, making the derivation of an adaptive robust control law unattainable. To address this coupling issue, we aim to separate the learning of uncertainty parameters from robust control, thereby leading to a tractable adaptive robust control law.  
\end{myRemark}

\section{Controller Design}\label{section3}
In this section, we describe the design of our adaptive robust control for an uncertain system disturbed by ellipsoid-set-bounded noises in detail. In Section \ref{section3_1}, we present the ellipsoid-set learning algorithm for updating the ellipsoid set of system states. In Section \ref{section3_2}, we integrate the updated state set in the control and provide a min-max robust control law with $N_p$-steps prediction for each candidate system parameter $\bm{\theta}_i$. We employ Bayes' theorem to update the \textit{a-posteriori} probability for each candidate parameter $\bm{\theta}_i$, where the \textit{a-posteriori} probability of the true parameter is likely to converge to $1$. We provide the convergence proof in Section \ref{section3_3}. In Section \ref{section3_4}, we derive the adaptive robust control law by integrating the robust control laws for each $\bm{\theta}_i$ weighted by their \textit{a-posteriori} probabilities. In Section \ref{section3_5}, we extend the proposed approach to solve the linear quadratic tracking problem. 

\subsection{The Ellipsoid-Set Learning for State Estimation} \label{section3_1}
In this section, we detail the ellipsoid-set learning algorithm for state estimation of systems whose process and observation noises are unknown, yet bounded in ellipsoid sets. We assume the ellipsoid set of initial state $E(\bar{\bm{x}}_0,\bm{P}_0)$ is \textit{a-priori} known, and this initial ellipsoid set is sufficiently large to cover all possible worst cases. Following the main technique in \cite{maksarov1996state}, we utilize the observation $\bm{y}(k)$ to learn the ellipsoid set for the system state in each control iteration. This ellipsoid-set learning leads to reduced state set size, resulting in the reduction of uncertainty in the control. In the following, we derive the learning algorithm for the ellipsoid set of states under the candidate parameter $\bm{\theta}_i$. 

Let the system with the parameter $\bm{\theta}_i$ be
\begin{equation}\label{subsys1}
	\begin{aligned}
		\bm{x}(k+1)=\bm{A}(k,\bm{\theta}_i)\bm{x}(k)+\bm{B}(k,\bm{\theta}_i)\bm{u}(k)+\bm{w}(k),
	\end{aligned}
\end{equation}
\begin{equation}\label{subsys2}
	\begin{aligned}
		\bm{y}(k) = \bm{C}(k,\bm{\theta}_i)\bm{x}(k)+\bm{v}(k),
	\end{aligned}
\end{equation}
where the parameter $\bm{\theta}_i$ belongs to the parameter space $\bm{\Theta}$ defined in \eqref{Theta_set}.
We define the updated ellipsoidal set for the state $\bm{x}(k)$ with the parameter $\bm{\theta}_i$ at the instant $k$ as
\begin{equation}\label{}
\begin{aligned}
    \mathcal{X}(k|k,\bm{\theta}_i) = E(\hat{\bm{x}}(k|k,\bm{\theta}_i),\bm{P}(k|k,\bm{\theta}_i))  ,
\end{aligned}
\end{equation}
where the ellipsoid center $\hat{\bm{x}}(k|k,\bm{\theta}_i)$ and shape matrix $\bm{P}(k|k,\bm{\theta}_i)$ are obtained by utilizing the observation $\bm{y}(k)$ to update the \textit{a-priori} ellipsoid set $\mathcal{X}(k-1|k-1,\bm{\theta}_i)=E(\hat{\bm{x}}(k-1|k-1,\bm{\theta}_i),\bm{P}(k-1|k-1,\bm{\theta}_i)) $.

\begin{figure}[]
\centering
\includegraphics[width=0.5\textwidth]{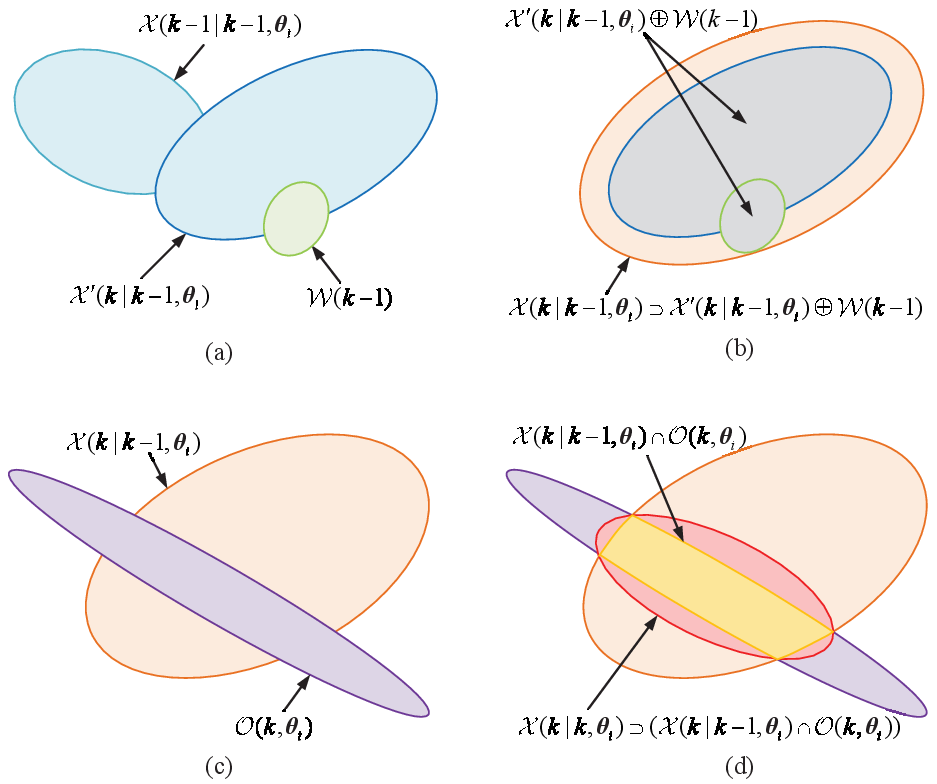}\\
\caption{Geometrical illustration of the learning progress of the ellipsoid set for the system state at the \textit{k}-th instant. }\label{Fig1}
\end{figure}

We present a geometrical illustration of the algorithm of ellipsoid-set learning in Fig. \ref{Fig1}, taking an example of a two-dimensional system state. In Fig. \ref{Fig1}(a), we linearly transform the ellipsoid set $\mathcal{X}(k-1|k-1,\bm{\theta}_i)$ to the set $\mathcal{X}'(k|k-1,\bm{\theta}_i) = E(\hat{\bm{x}}'(k|k-1,\bm{\theta}_i), \bm{P}'(k|k-1,\bm{\theta}_i))$, where the center and shape matrix are calculated respectively as
\begin{equation}\label{}
    \begin{aligned}
        \hat{\bm{x}}'(k|k-1,\bm{\theta}_i) = \bm{A}(k-1,\bm{\theta}_i)\hat{\bm{x}}(k-1|k-1,\bm{\theta}_i)\\
        +\bm{B}(k-1,\bm{\theta}_i)\bm{u}(k-1),
    \end{aligned}
\end{equation}
and
\begin{equation}\label{}
    \begin{aligned}
        \bm{P}'(k|k-1,\bm{\theta}_i) = \bm{A}(k-1,\bm{\theta}_i)\bm{P}(k-1|k-1,\bm{\theta}_i)\\
        \bm{A}^T(k-1,\bm{\theta}_i).
    \end{aligned}
\end{equation}
In Fig. \ref{Fig1}(b), we search a minimal ellipsoid $\mathcal{X}(k|k-1,\bm{\theta}_i)$ that can cover the vector sum of the set $\mathcal{X}'(k|k-1,\bm{\theta}_i)$ and process noise set $\mathcal{W}(k-1)=E(\bm{0},\bm{P}_{w}(k-1))$. We represent the ellipsoid $\mathcal{X}(k|k-1,\bm{\theta}_i)$ as
\begin{equation}\label{}
    \begin{aligned}
        \mathcal{X}(k|k-1,\bm{\theta}_i) \supset \mathcal{X}'(k|k-1,\bm{\theta}_i) \oplus \mathcal{W}(k-1),
    \end{aligned}
\end{equation}
where $\oplus$ is the Minkowski sum defined in our \textit{Notation}. We show the derivation of $ \mathcal{X}(k|k-1,\bm{\theta}_i)$ in Lemma \ref{lemma_sum} below. 

\begin{myLemma}\label{lemma_sum}
$\forall p(k,\bm{\theta}_i)\in(0,\infty)$, we have $\mathcal{X}(k|k-1,\bm{\theta}_i) \supset \mathcal{X}'(k|k-1,\bm{\theta}_i) \oplus \mathcal{W}(k-1)$, where the ellipsoid center and shape matrix of $\mathcal{X}(k|k-1,\bm{\theta}_i)=E(\hat{\bm{x}}(k|k-1,\bm{\theta}_i), \bm{P}(k|k-1,\bm{\theta}_i))$ are 
\begin{equation}\label{learning1}
	\begin{aligned}
		\hat{\bm{x}}(k|k-1,\bm{\theta}_i)=\bm{A}(k-1,\bm{\theta}_i)\hat{\bm{x}}(k-1|k-1,\bm{\theta}_i)\\
        +\bm{B}(k-1,\bm{\theta}_i)\bm{u}(k-1),
	\end{aligned}
\end{equation}
and 
\begin{equation}\label{learning2}
	\begin{aligned}
		\bm{P}(k|k-1,\bm{\theta}_i)=(p^{-1}(k,\bm{\theta}_i)+1)\bm{A}(k-1,\bm{\theta}_i)\\
		\bm{P}(k-1|k-1,\bm{\theta}_i)\bm{A}^T(k-1,\bm{\theta}_i)\\
            +(p(k,\bm{\theta}_i)+1)\bm{P}_w(k-1).
	\end{aligned}
\end{equation}
The optimal scalar parameter $p(k,\bm{\theta}_i)\in(0,\infty)$ that minimize the size of ellipsoid $E(\hat{\bm{x}}(k|k-1,\bm{\theta}_i), \bm{P}(k|k-1,\bm{\theta}_i))$ is 
\begin{equation}\label{p_i}
    p(k,\bm{\theta}_i) = \left[ \frac{\tr(\bm{P}'(k|k-1,\bm{\theta}_i))}{\tr(\bm{P}_w(k-1))} \right] ^{1/2}.
\end{equation}
\end{myLemma}
\begin{proof}
We detail the proof in Appendix \ref{Appendix_Lemma_1}.
\end{proof}

Fig. \ref{Fig1}(c) shows the intersection of the updated set $\mathcal{X}(k|k-1,\bm{\theta}_i)$ and the observation set
\begin{equation}\label{}
    \begin{aligned}
        \mathcal{O}(k)=\{\bm{x}(k)\in\mathbb{R}^n |[\bm{y}(k)-\bm{C}(k,\bm{\theta}_i)\bm{x}(k)]^T\bm{P}_{v}^{-1}(k)\\
        [\bm{y}(k)-\bm{C}(k,\bm{\theta}_i)\bm{x}(k)]\},
    \end{aligned}
\end{equation}
which can be rewritten as the standard form as $\mathcal{O}(k)=E(\bm{C}^{-1}_r(k,\bm{\theta}_i)\bm{y}(k),\bm{C}^{-1}_r(k,\bm{\theta}_i)\bm{P}_{v}(k)(\bm{C}^{-1}_r(k,\bm{\theta}_i))^T)$, and we provide the detailed derivation in Appendix \ref{Appendix_obs}.

Fig. \ref{Fig1}(d) shows the ellipsoid set update by a minimal ellipsoid set $\mathcal{X}(k|k,\bm{\theta}_i)$ that covers the overlap between $\mathcal{X}(k|k-1,\bm{\theta}_i)$ and $\mathcal{O}(k)$. We represent the process where the ellipsoid set is updated using the observation set $\mathcal{O}(k)$ as 
\begin{equation}\label{}
	\begin{aligned}
		\mathcal{X}(k|k,\bm{\theta}_i) \supset \mathcal{X}(k|k-1,\bm{\theta}_i) \cap \mathcal{O}(k).
	\end{aligned}
\end{equation}
We show the derivation of $ \mathcal{X}(k|k,\bm{\theta}_i)$ in Lemma \ref{lemma_intersect} below. 

\begin{myLemma}\label{lemma_intersect}
$\forall q(k,\bm{\theta}_i) \in [0,\infty)$, we have $\mathcal{X}(k|k,\bm{\theta}_i) \supset \mathcal{X}(k|k-1,\bm{\theta}_i) \cap \mathcal{O}(k)$ with the ellipsoid center and shape matrix of $\mathcal{X}(k|k,\bm{\theta}_i)=E(\hat{\bm{x}}(k|k,\bm{\theta}_i), \bm{P}(k|k,\bm{\theta}_i))$ being chosen as 
\begin{equation}\label{learning4}
	\begin{aligned}
		\hat{\bm{x}}(k|k,\bm{\theta}_i)=\hat{\bm{x}}(k|k-1,\bm{\theta}_i)+\bm{K}(k,\bm{\theta}_i)\bm{\epsilon}(k,\bm{\theta}_i)
	\end{aligned}
\end{equation}
and 
\begin{equation}\label{learning5}
	\begin{aligned}
		\bm{P}(k|k,\bm{\theta}_i)=\beta(k,\bm{\theta}_i)\{ \left[ \bm{I}-\bm{K}(k,\bm{\theta}_i)\bm{C}(k,\bm{\theta}_i) \right]\\
            \bm{P}(k|k-1,\bm{\theta}_i) \left[ \bm{I}-\bm{K}(k,\bm{\theta}_i)\bm{C}(k,\bm{\theta}_i) \right]^T \\
            + q^{-1}(k,\bm{\theta}_i)\bm{K}(k,\bm{\theta}_i)\bm{P}_v(k)
        \bm{K}^T(k,\bm{\theta}_i) \},
	\end{aligned}
\end{equation}
where the gain $\bm{K}(k,\bm{\theta}_i)\in \mathbb{R}^{n \times m}$, innovation $\bm{\epsilon}(k,\bm{\theta}_i)\in \mathbb{R}^{m} $ and parameter $\beta(k,\bm{\theta}_i) \in \mathbb{R}$ are respectively defined as
\begin{equation}\label{learning3}
	\begin{aligned}
		\bm{K}(k,\bm{\theta}_i)=\bm{P}(k|k-1,\bm{\theta}_i)\bm{C}^T(k,\bm{\theta}_i)
		\left[ \bm{C}(k,\bm{\theta}_i) \right. \\
       \left. \bm{P}(k|k-1,\bm{\theta}_i)\bm{C}^T(k,\bm{\theta}_i)  +q^{-1}(k,\bm{\theta}_i)\bm{P}_v(k) \right] ^{-1},
	\end{aligned}
\end{equation}
\begin{equation}\label{learning7}
    \begin{aligned}
        \bm{\epsilon}(k,\bm{\theta}_i) = \bm{y}(k)-\bm{C}(k,\bm{\theta}_i)\hat{\bm{x}}(k|k-1,\bm{\theta}_i),
    \end{aligned}
\end{equation}
\begin{equation}\label{learning6}
	\begin{aligned}
		\beta(k,\bm{\theta}_i)=1+q(k,\bm{\theta}_i)-\bm{\epsilon}^T(k,\bm{\theta}_i)
		\left[ q^{-1}(k,\bm{\theta}_i)\bm{P}_v(k)\right.  \\
		\left. +\bm{C}(k,\bm{\theta}_i)\bm{P}(k|k-1,\bm{\theta}_i)\bm{C}^T(k,\bm{\theta}_i) \right]^{-1}\bm{\epsilon}(k,\bm{\theta}_i).
	\end{aligned}
\end{equation}
The optimal scalar parameter $q(k,\bm{\theta}_i)$ is obtained by solving:
\begin{equation}\label{solveq}
	\begin{aligned}
		\frac{\sum_{j=1}^{n} d_j(k,\bm{\theta}_i) \frac{\lambda_j(k,\bm{\theta}_i)}{(1+q(k,\bm{\theta}_i)\lambda_j(k,\bm{\theta}_i))^2} }{\sum_{j=1}^{n} d_j(k,\bm{\theta}_i) \frac{1}{1+q(k,\bm{\theta}_i)\lambda_j(k,\bm{\theta}_i)}  }  = \frac{\beta'(k,\bm{\theta}_i)}{\beta(k,\bm{\theta}_i)},
	\end{aligned}
\end{equation}
where $d_j(k,\bm{\theta}_i)$ is the $j$-th diagonal element of matrix $\bm{V}^{-1}\bm{P}(k|k-1,\bm{\theta}_i)\bm{V}$ with $\bm{V}\bm{D}\bm{V}^{-1}$ being the diagonal decomposition of $\bm{P}(k|k-1,\bm{\theta}_i)\bm{C}(k,\bm{\theta}_i)^T\bm{P}_v^{-1}(k)\bm{C}(k,\bm{\theta}_i)$. $\bm{D}$ is the diagonal matrix containing eigenvalues of matrix $\bm{P}(k|k-1,\bm{\theta}_i)\bm{C}(k,\bm{\theta}_i)^T\bm{P}_v^{-1}(k)\bm{C}(k,\bm{\theta}_i)$, where $\lambda_j(k,\bm{\theta}_i)$ is the $j$-th eigenvalue, and $n$ is the number of the eigenvalues. $\beta'(k,\bm{\theta}_i)$ is the partial derivative of $\beta(k,\bm{\theta}_i)$ with respect to $q(k,\bm{\theta}_i)$, calculated by
\begin{equation}
	\begin{aligned}
		\beta'(k,\bm{\theta}_i)=1-\bm{\epsilon}^T(k,\bm{\theta}_i)\left[ q^{-1}(k,\bm{\theta}_i)\bm{P}_v(k)+\bm{C}(k,\bm{\theta}_i) \right.\\
        \left. \bm{P}(k|k-1,\bm{\theta}_i)\bm{C}^T(k,\bm{\theta}_i) \right]^{-1}q^{-2}(k,\bm{\theta}_i)\bm{P}_v(k)\left[ q^{-1}(k,\bm{\theta}_i) \right.\\
        \left. \bm{P}_v(k)+\bm{C}(k,\bm{\theta}_i)\bm{P}(k|k-1,\bm{\theta}_i)\bm{C}^T(k,\bm{\theta}_i) \right]^{-1} \bm{\epsilon}(k,\bm{\theta}_i),
	\end{aligned}
\end{equation}
Note that $q(k,\bm{\theta}_i)=0$, if
\begin{equation}\label{cond_nosol}
	\begin{aligned}
		\left( 1-\bm{\epsilon}^T(k,\bm{\theta}_i)\bm{P}_v^{-1}(k)\bm{\epsilon}(k,\bm{\theta}_i)\right)\tr\left( \bm{P}(k|k-1,\bm{\theta}_i)\right) \\
		 -\tr\left( \bm{P}(k|k-1,\bm{\theta}_i)\bm{C}^T(k,\bm{\theta}_i)\bm{P}_v^{-1}(k)\bm{C}(k,\bm{\theta}_i) \right.\\
        \left. \bm{P}(k|k-1,\bm{\theta}_i) \right)  >0.
	\end{aligned}
\end{equation}
\end{myLemma}
\begin{proof}
We detail the proof in Appendix \ref{Appendix_Lemma_2}.
\end{proof}

\begin{myRemark}
The ellipsoid-set learning algorithm exhibits structural similarities to Kalman filter. In particular, both methods predict and update the state uncertainty (described by Gaussian distributions in Kalman filter, and by ellipsoidal sets in ellipsoid-set learning) upon the state evolution equation and real-time observations. In Kalman filter, the covariance matrix of the state estimation corresponds to a certain confidence region of the Gaussian distribution, which can be intuitively visualized as an ellipsoid. However, the covariance matrix updated by Kalman filter provides an over-optimistic ellipsoid when it comes to bounded noise. In contrast, the ellipsoid-set learning algorithm estimates the set that covers all possible system states, which is more robust than the classical Kalman filter. Thus, the ellipsoid-set learning algorithm is suitable for control applications where the system states are required to be confined within boundaries, such as the robust control for safety-critical systems. 
\end{myRemark}

\subsection{Robust N$_{p}$-Steps Predictive Control}\label{section3_2}
This section presents a robust control that integrates the ellipsoid-set learning algorithm as shown in Section \ref{section3_1}. The control in this section avoids open-loop control policies derived based on a fixed initial state set~\cite{bertsimas2007constrained}. Instead, we update the state set iteratively and use it to design a closed-loop control policy. Thus we reduce the conservativeness of the control with the state set learned and narrowed through the ellipsoid-set learning. Below we describe the derivation of our robust control in detail. 

According to \textit{Assumption} \ref{Assump_obser}, the true value of the state $\bm{x}(k)$ is not directly accessible. We utilize the observation $\bm{y}(k)$ to estimate the ellipsoid set of the state by applying the ellipsoid-set learning algorithm described in Section \ref{section3_1}. The true value of the state is confined in the updated ellipsoid set $E(\hat{\bm{x}}(k|k,\bm{\theta}_i),\bm{P}(k|k,\bm{\theta}_i))$, where $\hat{\bm{x}}(k|k,\bm{\theta}_i)$ is center and $\bm{P}(k|k,\bm{\theta}_i)$ is the shape matrix of the ellipsoid. We add the updated set as a constraint for the state $\bm{x}(k)$ in formulating the problem $(\mathcal{P}_i)$, and solve it to derive the closed-loop control law. The problem $(\mathcal{P}_i)$ takes the form: 
\begin{equation}\label{prob_pred}
    \begin{aligned}
        (\mathcal{P}_i): \quad  &\min_{\{{u}(t,\bm{\theta}_i)\} \vert_{t=k}^{k+N_p-1}} \max_{ \{\bm{w}(k), \bm{v}(k),\bm{x}(k)\} } J_k \\
        \quad s.t. \quad &\bm{x}(k+1)=\bm{A}(k,\bm{\theta}_i)\bm{x}(k)+\bm{B}(k,\bm{\theta}_i)\bm{u}(k)+\bm{w}(k),\\
        \quad \quad \quad \quad & \bm{y}(k) = \bm{C}(k,\bm{\theta}_i)\bm{x}(k)+\bm{v}(k), \\
        \quad \quad \quad \quad & \bm{w}(k) \in E(\bm{0},\bm{P}_w(k)), \\
        \quad \quad \quad \quad & \bm{v}(k) \in E(\bm{0},\bm{P}_v(k)), \\
        \quad \quad \quad \quad & \bm{x}(k)\in E(\hat{\bm{x}}(k|k,\bm{\theta}_i),\bm{P}(k|k,\bm{\theta}_i)).
    \end{aligned}
\end{equation}
Note that rather than optimizing the cost function within the whole horizon $N$ as \eqref{prob1}, we reduce the optimization horizon in solving the problem $(\mathcal{P}_i)$ to reduce the computational complexity. Particularly, we apply the receding horizon approximation (also called model predictive control) to reduce the whole horizon optimization problem to the problem with $N_p$-step prediction and optimization $(\mathcal{P}_i)$~\cite{schwenzer2021review}.
The cost function with the prediction horizon $[k,k+N_p-1]$ is 
\begin{equation}\label{J_N}
	\begin{aligned}
		J_k = \sum_{t= k}^{k+N_p-1} \left[\bm{x}^T(t)\bm{Q}(t)\bm{x}(t) +  \bm{u}^T(t)\bm{R}(t)\bm{u}(t) \right].
	\end{aligned}
\end{equation}
Note that the control goal is finding an optimal control sequence $\{{u}(t,\bm{\theta}_i)\} \vert_{t=k}^{k+N_p-1}$ that minimizes the above cost function at each instant $k$. 

In the following, we show how we solve the optimization problem $(\mathcal{P}_i)$. Lemma~\ref{lemma3} rewrites the cost function \eqref{J_N} as a quadratic function of the control sequence. Due to that the true value of state $\bm{x}(k)$ is not attainable, we introduce the estimated ellipsoid set of state $E(\hat{\bm{x}}(k|k,\bm{\theta}_i),\bm{P}(k|k,\bm{\theta}_i))$ into the cost function \eqref{J_N} in Lemma \ref{lemma4}. In Theorem~\ref{theorem1}, we formulate the optimization problem $(\mathcal{P}_i)$ into a solvable SDP problem. 

\begin{myLemma}\label{lemma3}
The cost function $J_k$ defined in (\ref{J_N}) can be expressed as
\begin{equation}\label{J_t_x}
    \begin{aligned}
        &J_k =  \bm{U}^T_k\mathscr{B}_k(\bm{\theta}_i)\bm{U}_k+2\bm{b}^T_k(\bm{\theta}_i)\bm{U}_k
     +\bm{W}^T_k\mathscr{C}_k(\bm{\theta}_i)\bm{W}_k \\ &+2\bm{c}^T_k(\bm{\theta}_i) \bm{W}_k +2\bm{U}^T_k\mathscr{D}_k(\bm{\theta}_i)\bm{W}_k +\bm{x}^T(k)\mathscr{A}_k(\bm{\theta}_i)\bm{x}(k) ,
    \end{aligned}
\end{equation}
where the control vector $\bm{U}_k \in \mathbb{R}^{N_p\cdot r}$ and process noise vector $\bm{W}_k \in \mathbb{R}^{N_p\cdot n } $ are respectively defined as 
\begin{equation}\label{U_vec}
    \begin{aligned}
        \bm{U}_k= [ \bm{u}^T(k) \quad \bm{u}^T(k+1) \quad \cdots \quad \bm{u}^T(k+N_p-1]^T ,
    \end{aligned}
\end{equation}
\begin{equation}\label{W_vec}
        \bm{W}_k= [  \bm{w}^T(k) \quad \bm{w}^T(k+1) \quad \cdots  \quad \bm{w}^T(k+{\tiny }N_p-1)  ]^T.
\end{equation}
The parameter matrices $\mathscr{A}_k(\bm{\theta}_i) \in \mathbb{R}^{n\times n} $, $\mathscr{B}_k(\bm{\theta}_i) \in \mathbb{R}^{N_p\cdot r \times N_p\cdot r} $, $\mathscr{C}_k(\bm{\theta}_i) \in \mathbb{R}^{N_p\cdot n \times N_p\cdot n} $, and $\mathscr{D}_k(\bm{\theta}_i) \in \mathbb{R}^{N_p\cdot r \times N_p\cdot n} $ are respectively defined as
\begin{equation}\label{A_math}
    \begin{aligned}
        \mathscr{A}_k(\bm{\theta}_i) = \sum_{t=k}^{k+N_p-1}(\tilde{\bm{A}}^{t-1}_k)^T_i\bm{Q}(t)(\tilde{\bm{A}}^{t-1}_k)_i,
    \end{aligned}
\end{equation}
\begin{equation}\label{B_math}
    \begin{aligned}
        &\mathscr{B}_k(\bm{\theta}_i) = \sum_{t=k+1}^{k+N_p-1}(\tilde{\bm{B}}^{t-1}_k)^T_i \bm{Q}(t) (\tilde{\bm{B}}^{t-1}_k)_i \\ &+\diag(\bm{R}(k),\bm{R}(k+1), \cdots, \bm{R}(k+N_p-1)),
    \end{aligned}
\end{equation}
\begin{equation}\label{C_math}
    \begin{aligned}
        \mathscr{C}_k(\bm{\theta}_i) = \sum_{t=k+1}^{k+N_p-1}(\tilde{\bm{C}}^{t-1}_k)^T_i \bm{Q}(t) (\tilde{\bm{C}}^{t-1}_k)_i,
    \end{aligned}
\end{equation}
\begin{equation}\label{D_math}
    \begin{aligned}
        \mathscr{D}_k(\bm{\theta}_i) = \sum_{t=k+1}^{k+N_p-1} (\tilde{\bm{B}}^{t-1}_k)^T_i \bm{Q}(t) (\tilde{\bm{C}}^{t-1}_k)_i.
    \end{aligned}
\end{equation}
The parameter vectors $\bm{b}_k(\bm{\theta}_i) \in \mathbb{R}^{N_p\cdot r} $ and $\bm{c}_k(\bm{\theta}_i) \in \mathbb{R}^{N_p\cdot n}$ are respectively defined as
\begin{equation}\label{b_math}
    \bm{b}_k(\bm{\theta}_i) = \sum_{t=k+1}^{k+N_p-1}(\tilde{\bm{B}}^{t-1}_k)^T_i\bm{Q}(t)(\tilde{\bm{A}}^{t-1}_k)_i\bm{x}(k),
\end{equation}
\begin{equation}\label{c_math}
    \begin{aligned}
        \bm{c}_k(\bm{\theta}_i) = \sum_{t=k+1}^{k+N_p-1}(\tilde{\bm{C}}^{t-1}_k)^T_i\bm{Q}(t)(\tilde{\bm{A}}^{t-1}_k)_i\bm{x}(k).
    \end{aligned}
\end{equation}
The matrices $(\tilde{\bm{A}}^{t-1}_k)_i \in \mathbb{R}^{n\times n} $, $(\tilde{\bm{B}}^{t-1}_k)_i \in \mathbb{R}^{n \times N_p\cdot r} $, and $(\tilde{\bm{C}}^{t-1}_k)_i \in \mathbb{R}^{n\times N_p\cdot n } $ are calculated by
\begin{equation}\label{A_wave}
    \begin{aligned}
        (\tilde{\bm{A}}^{t-1}_k)_i=\left\lbrace \begin{array}{rcl}
            \prod_{l=k}^{t-1}\bm{A}(l,\bm{\theta}_i) & \mbox{for} &  t \geqslant k+1  \\
            1& \mbox{for} & t=k
        \end{array} \right.,
    \end{aligned}
\end{equation}
\begin{equation}\label{B_wave}
    \begin{aligned}
        (\tilde{\bm{B}}^{t-1}_k)_i= [(\tilde{\bm{A}}^{t-1}_{k+1})_i\bm{B}(k,\bm{\theta}_i) \quad (\tilde{\bm{A}}^{t-1}_{k+2})_i\bm{B}(k+1,\bm{\theta}_i)  \\
        \quad \cdots \quad \bm{B}(t-1,\bm{\theta}_i) \quad \bm{0}_{n \times (N_p+k-t)\cdot r} ] ,
    \end{aligned}
\end{equation}
\begin{equation}\label{C_wave}
        (\tilde{\bm{C}}^{t-1}_k)_i = [ (\tilde{\bm{A}}^{t-1}_{k+1})_i \quad (\tilde{\bm{A}}^{t-1}_{k+2})_i \quad \cdots   \quad  \bm{0}_{n \times (N_p+k-t)\cdot n}],
\end{equation}
and $t = k,k+1,\cdots , k+N_p-1$.
\end{myLemma}
\begin{proof}
The proof can be found in our previous work \cite{ma2024robust}.
\end{proof}

According to Lemma~\ref{lemma_sum} and \ref{lemma_intersect}, we obtain the ellipsoid set of the state at every instant $k$
\begin{equation}\label{stateset}
	\begin{aligned}
		\mathcal{X}(k|k,\bm{\theta}_i)&= E(\hat{\bm{x}}(k|k,\bm{\theta}_i),\bm{P}(k|k,\bm{\theta}_i))\\
  &=\{\bm{x}\in \mathbb{R}^n|(\bm{x}(k)-\hat{\bm{x}}(k|k,\bm{\theta}_i))^T\\
  &\bm{P}^{-1} (k|k,\bm{\theta}_i)(\bm{x}(k)-\hat{\bm{x}}(k|k,\bm{\theta}_i))\leq 1\}.
	\end{aligned}
\end{equation}
Define the estimation error of the state at the instant $k$ as 
\begin{equation}\label{eta_def}
	\begin{aligned}
		\bm{\eta}(k,\bm{\theta}_i) = \bm{x}(k)-\hat{\bm{x}}(k|k,\bm{\theta}_i).
	\end{aligned}
\end{equation}
Substituting (\ref{eta_def}) into (\ref{stateset}), we have the estimation error confined in the following ellipsoid set
\begin{equation}\label{eta_set}
	\begin{aligned}
		\bm{\eta}(k,\bm{\theta}_i) \in E(\bm{0},\bm{P}(k|k,\bm{\theta}_i)).
	\end{aligned}
\end{equation}
In Lemma~\ref{lemma4} below, we rewrite the quadratic form of the cost function $J_k$ by introducing the disturbance vector $\bm{\xi}_k^T(\bm{\theta}_i) = \left[ \begin{array}{cc} \bm{\eta}^T(k,\bm{\theta}_i) & \bm{W}_k^T \end{array} \right]$, where the estimation error is added. This facilitates us to apply the S-procedure and Schur complement \cite{boyd2004convex}, to reformulate the problem $(\mathcal{P}_i)$ to an SDP problem.

\begin{myLemma}\label{lemma4}
With the state estimation error given by (\ref{eta_def}), we can rewrite the cost function $J_k$ in (\ref{J_t_x}) as
\begin{equation}\label{J_t_xhat}
    \begin{aligned}
        J_k = &\bm{U}_k^T\mathscr{B}_k(\bm{\theta}_i)\bm{U}_k+2\hat{\bm{b}}_k^T(\bm{\theta}_i)\bm{U}_k +\bm{\xi}_k^T(\bm{\theta}_i)\hat{\mathscr{C}}_k(\bm{\theta}_i)\bm{\xi}_k(\bm{\theta}_i)\\
        &+2\hat{\bm{c}}_k^T(\bm{\theta}_i) \bm{\xi}_k(\bm{\theta}_i)+2\bm{U}_k^T\hat{\mathscr{D}}_k(\bm{\theta}_i)\bm{\xi}_k(\bm{\theta}_i)\\
        &+\hat{\bm{x}}^T(k|k,\bm{\theta}_i)\mathscr{A}_k(\bm{\theta}_i) \hat{\bm{x}}(k|k,\bm{\theta}_i)
    \end{aligned}
\end{equation}
where $\hat{\bm{b}}_k(\bm{\theta}_i) \in \mathbb{R}^{N_p\cdot r } $, $\hat{\mathscr{C}}_k(\bm{\theta}_i) \in \mathbb{R}^{(N_p+1)\cdot n \times (N_p+1)\cdot n}$, $\hat{\bm{c}}_k(\bm{\theta}_i) \in \mathbb{R}^{(N_p+1)\cdot n } $, $\hat{\mathscr{D}}_k(\bm{\theta}_i) \in \mathbb{R}^{N_p\cdot r \times (N_p+1)\cdot n} $, $\bm{\xi}_k \in \mathbb{R}^{(N_p+1)\cdot n}$ are, respectively, defined as
\begin{equation}
        \hat{\bm{b}}_k(\bm{\theta}_i) = \bm{S}_{B}(\bm{\theta}_i) \hat{\bm{x}}(k|k),
\end{equation}
\begin{equation}
    \begin{aligned}
        \hat{\mathscr{C}}_k(\bm{\theta}_i)= \left[ \begin{array}{cc} \mathscr{A}_k(\bm{\theta}_i) & \bm{S}_{C}^T(\bm{\theta}_i) \\ \bm{S}_{C}(\bm{\theta}_i) & \mathscr{C}_k(\bm{\theta}_i) \end{array} \right] , 
    \end{aligned}
\end{equation}
\begin{equation}
    \begin{aligned}
        \hat{\bm{c}}_k(\bm{\theta}_i)=\left[ \begin{array}{c} \mathscr{A}_k(\bm{\theta}_i)\hat{\bm{x}}(k|k) \\ \bm{S}_{C}(\bm{\theta}_i) \hat{\bm{x}}(k|k) \end{array} \right] ,
    \end{aligned}
\end{equation}
\begin{equation}
    \begin{aligned}
        \hat{\mathscr{D}}_k(\bm{\theta}_i)=
\left[ \bm{S}_{B}(\bm{\theta}_i) \quad \mathscr{D}_k(\bm{\theta}_i) \right] ,
    \end{aligned}
\end{equation}
\begin{equation}
    \begin{aligned}
        \bm{\xi}_k^T(\bm{\theta}_i) = \left[ \begin{array}{cc} \bm{\eta}^T(k,\bm{\theta}_i) & \bm{W}_k^T \end{array} \right] .
    \end{aligned}
\end{equation}	
\end{myLemma}
\begin{proof}
The proof can be found in our previous work \cite{ma2024robust}.
\end{proof}

\begin{myTheo}\label{theorem1}
The problem $(\mathcal{P}_i)$ in (\ref{prob_pred}) is equivalent to the following formulation:
\begin{equation}\label{LMI_Theo}
    \begin{aligned}
        &(\mathcal{P}_i): \quad  \min \quad \rho_k  \\
        \quad s.t.  &
        \left[ \begin{array}{ccc} \bm{I} & \bm{z}_k & \bm{F}_k(\bm{\theta}_i) \\ \bm{z}_k^T & \rho_k-\tau_1-\tau_2 N_p & -\bm{h}_k^T(\bm{\theta}_i) \\ \bm{F}_k^T(\bm{\theta}_i) &  -\bm{h}_k(\bm{\theta}_i) & \bm{G}_k(\bm{\theta}_i) \end{array} \right]   \geq 0,
    \end{aligned}
\end{equation}
where $\tau_1 \geq 0$, $\tau_2 \geq 0$, $\bm{z}_k\in \mathbb{R}^{n_p\cdot r}$, and $\rho_k \in \mathbb{R}$ are decision variables. The vector $\bm{h}_k(\bm{\theta}_i)$ and the matrix $\bm{F}_k(\bm{\theta}_i)$, $\bm{M^{\eta}}(k,\bm{\theta}_i)$, $\bm{M}_k^{W}$, $\bm{G}_k(\bm{\theta}_i)$ are respectively defined as
\begin{equation}\label{h_vector}
    \begin{aligned}
        \bm{h}_k(\bm{\theta}_i) = \hat{\bm{c}}_k(\bm{\theta}_i) - \hat{\mathscr{D}}_k^T(\bm{\theta}_i)\mathscr{B}_k^{-1}(\bm{\theta}_i)\hat{\bm{b}}_k(\bm{\theta}_i),
    \end{aligned}
\end{equation}
\begin{equation}\label{F_matrix}
    \begin{aligned}
        \bm{F}_k(\bm{\theta}_i) = \mathscr{B}_k^{-1/2}(\bm{\theta}_i)\hat{\mathscr{D}}_k(\bm{\theta}_i).
    \end{aligned}
\end{equation}
\begin{equation}\label{Matrix_eta}
    \begin{aligned}
        \bm{M^{\eta}}(k,\bm{\theta}_i) = \left[ \begin{array}{cc} \bm{P}^{-1}(k,\bm{\theta}_i) & 0 \\ 0 & 0 \end{array} \right] ,
    \end{aligned}
\end{equation}
\begin{equation}\label{Matrix_W}
    \begin{aligned}
        \bm{M}_k^{W} = \left[ \begin{array}{cc} 0 & 0 \\ 0 & \bm{P}_k^W \end{array} \right] ,
    \end{aligned}
\end{equation}
\begin{equation}
    \begin{aligned}
        \bm{G}_k(\bm{\theta}_i) = & -\hat{\mathscr{C}}_k(\bm{\theta}_i)+\tau_1\bm{M^{\eta}}(k,\bm{\theta}_i)+\tau_2\bm{M}_k^{W}\\
        &+\bm{F}_t^T(\bm{\theta}_i)\bm{F}_t(\bm{\theta}_i) ,
    \end{aligned}
\end{equation}
where $\bm{P}_k^W = \diag(\bm{P}_w^{-1}(k), \cdots, \bm{P}_w^{-1}(k+N_p-1))$. We calculate the corresponding control vector defined in (\ref{U_vec}) by
\begin{equation}\label{U_vecotr}
    \begin{aligned}
        \bm{U}_k = \mathscr{B}_k^{-1/2}(\bm{\theta}_i)\bm{z}_k-\mathscr{B}_k^{-1}(\bm{\theta}_i)\hat{\bm{b}}_k(\bm{\theta}_i),
    \end{aligned}
\end{equation}
and the control law at the instant $k$ is the first element of $\bm{U}_k$:
\begin{equation}\label{equation_MPC_control_law}
    \begin{aligned}
        \bm{u}(k,\bm{\theta}_i) = \bm{U}_k(0).
    \end{aligned}
\end{equation}
\end{myTheo}
\begin{proof}
We provide the detailed proof in Appendix \ref{Appendix_Theo2}.
\end{proof}
Note that after the reformulation, now the problem $(\mathcal{P}_i)$ in Theorem~\ref{theorem1} is a typical SDP problem that can be solved efficiently by several available algorithms, like interior point methods~\cite{alizadeh1995interior}.

\subsection{Convergence of the Weights for Parameter Vectors }\label{section3_3}
The control law derivation in Section~\ref{section3_2} is based on the candidate parameter $\bm{\theta}_i$. As we describe in problem $(\mathcal{P})$, the uncertain parameter vector $\bm{\theta}$ takes value in the finite candidate set $\bm{\Theta}=\{\bm{\theta}_1, \bm{\theta}_2, \cdots, \bm{\theta}_s\}$, and we need to learn which admissible parameter vector $\bm{\theta}_i ( i=1,2,\cdots, s) $ is the true one for the control law derivation. To this end, we assign every $\bm{\theta}_i$ in $\bm{\Theta}$ with a corresponding weight, which can be calculated by the \textit{a-posteriori} probabilities. We adopt the Bayes' approach~\cite{DESHPANDE1973107} and update the \textit{a-posteriori} probability of the candidate parameters $\bm{\theta}_i$ using the observations at each instant. Below we present the learning algorithm for the weight of candidate parameter $\bm{\theta}_i$ and the convergence proof for the weight learning.

We define the weight for the system with parameter $\bm{\theta}_i$ at the $k$-th instant as $\pi(k,\bm{\theta}_i)$, and apply the \textit{a-posteriori} probability to calculate the weight, that is
\begin{equation}
    \begin{aligned}
        \pi(k,\bm{\theta}_i) = P(\bm{\theta}=\bm{\theta}_i \left| \bm{\mathfrak{I}}_k \right. ),\quad
         i = 1,2,\cdots, s, 
    \end{aligned}
\end{equation}
where $\bm{\mathfrak{I}}_k$ is the information set defined as
\begin{equation}
    \begin{aligned}
        \bm{\mathfrak{I}}_k = \{ \bm{u}(0), \cdots, \bm{u}(k-1), \bm{y}(1), \cdots, \bm{y}(k), \bm{\mathfrak{I}}_0\},
    \end{aligned}
\end{equation}
and $P(\bm{\theta}=\bm{\theta}_i \left| \bm{\mathfrak{I}}_k \right. )$ is the probability of the parameter vector $\bm{\theta}=\bm{\theta}_i$ given the evidence $\bm{\mathfrak{I}}_k$. According to Bayes' theorem, the weight, or the \textit{a-posteriori} probability of $\bm{\theta}_i$, can be updated by the following recursive equation 
\begin{equation}\label{PosProb}
	\begin{aligned}
		\pi(k,\bm{\theta}_i) = \frac{L(k,\bm{\theta}_i)}{\sum_{j=1}^{s}\pi(k-1,\bm{\theta}_j)L(k,\bm{\theta}_j)} \pi(k-1,\bm{\theta}_i)
	\end{aligned}
\end{equation}
where $L(k,\bm{\theta}_i)$ is the likelihood function. 
According to \textit{Assumption} \ref{Assump_noise} and \eqref{stateset}, only ellipsoid boundaries of states and noises are available, while no information about their exact probability distributions is required. Hence, we utilize some typical probability density functions to approximate and calculate the likelihood required in updating the weights in \eqref{PosProb}. For example, assuming the state $\bm{x}(k)$ and noise $\bm{v}(k)$ are uniformly distributed, we obtain that the observation $\bm{y}(k)$ is uniformly distributed in the ellipsoid set $\mathcal{Y}(k|k-1,\bm{\theta}_i)= E(\hat{\bm{y}}(k|k-1,\bm{\theta}_i), \bm{P}_y(k|k-1,\bm{\theta}_i))$, where 
\begin{equation}\label{y_hat}
    \hat{\bm{y}}(k|k-1,\bm{\theta}_i) = \bm{C}(k,\bm{\theta}_i)\hat{\bm{x}}(k|k-1,\bm{\theta}_i),
\end{equation}
\begin{equation}\label{P_y}
    \begin{aligned}
    \bm{P}_y(k|k-1,\bm{\theta}_i) = (p^{-1}_y(k,\bm{\theta}_i)+1)\bm{C}(k,\bm{\theta}_i)\\
    \bm{P}(k|k-1,\bm{\theta}_i)\bm{C}^T(k,\bm{\theta}_i)
    +(p_y(k,\bm{\theta}_i)+1)\bm{P}_v(k).
    \end{aligned}
\end{equation}
The volume of ellipsoid set $\mathcal{Y}(k|k-1,\bm{\theta}_i)$ is $\frac{\pi^{n/2}}{\Gamma(n/2+1)}\sqrt{\det( \bm{P}_y(k|k-1,\bm{\theta}_i) )}$. The likelihood function $L(k,\bm{\theta}_i)$ in (\ref{PosProb}) is the probability of $\bm{y}(k)$ given $\bm{\theta}_i$ and $\bm{\mathfrak{I}}_{k-1}$, which can be calculated by the inverse of the volume of ellipsoid $\mathcal{Y}(k|k-1,\bm{\theta}_i)$. Since the term $\frac{\Gamma(n/2+1)}{\pi^{n/2}}$ in $L(k,\bm{\theta}_i)$ and $L(k,\bm{\theta}_j)$ is the same, the likelihood $L(k,\bm{\theta}_i)$ can be simplified as
\begin{equation}\label{def_L}
	\begin{aligned}
		L(k,\bm{\theta}_i) = \left\{ \begin{array}{rl}  \frac{1}{\sqrt{\det( \bm{P}_y(k|k-1,\bm{\theta}_i) )}}, &  \bm{y}(k) \in \mathcal{Y}(k|k-1,\bm{\theta}_i) \\ 0, & \mbox{others} \end{array} \right.
	\end{aligned}
\end{equation}

We define the initial weights (\textit{a-prior} probabilities) of parameter vector $\bm{\theta}_i$ as
\begin{equation}\label{prob_init}
	\begin{aligned}
		\pi(0,\bm{\theta}_i) = P(\bm{\theta} = \bm{\theta}_i\left| \bm{\mathfrak{I}}_0 \right. ) , \quad i = 1, 2, \cdots , s,
	\end{aligned}
\end{equation}
and the sum of weights satisfy $\sum_{i=1}^{s}\pi(0,\bm{\theta}_i) = 1$.
Since the \textit{a-prior} probability is unknown in practice, we set the initial weights as the same value $1/s$.  With the initial weights, we can update the weights over time according to (\ref{PosProb}). As the part $\sum_{j=1}^{s}\pi(k-1,\bm{\theta}_j)L(k,\bm{\theta}_j)$ in (\ref{prob_init}) is the same for each $\pi(k,\bm{\theta}_i)$, we replace it by a constant as below
\begin{equation}\label{PosProb_a}
    \begin{aligned}
        \pi(k,\bm{\theta}_i) = \alpha(k) L(k,\bm{\theta}_i) \pi(k-1,\bm{\theta}_i)
    \end{aligned}
\end{equation}
where the constant $\alpha(k)$ is independent of $\bm{\theta}_i$, and the probabilities satisfy $\sum_{i=1}^{s}\pi(k,\bm{\theta}_i) = 1$.  
The posterior probability of $\bm{\theta}_i$ will converge to 1 if $\bm{\theta}_i$ is the true value of $\bm{\theta}$. We present the convergence proof for this weight updating algorithm in Theorem~\ref{thm2}.

\begin{myTheo}\label{thm2}
Suppose $\bm{\theta}_i$ is the true value of the unknown parameter vector $\bm{\theta}$ to be identified. We have that the posterior probability $\pi(k,\bm{\theta}_i) \rightarrow 1$ when $k\rightarrow \infty$, while the posterior probability $\pi(k,\bm{\theta}_j) \rightarrow 0$ when $k\rightarrow \infty$ for $\forall j \neq i $.	
\end{myTheo}

%\begin{proof}
%We provide the detailed proof in Appendix \ref{Appendix_Theo2}.
%\end{proof}

\begin{proof}
Combining \textit{Assumption} \ref{Assump_noise} and equation~\eqref{stateset}, we obtain that the observation $\bm{y}(k)$ is confined in the ellipsoid set $\mathcal{Y}(k|k-1,\bm{\theta}_i)= E(\hat{\bm{y}}(k|k-1,\bm{\theta}_i), \bm{P}_y(k|k-1,\bm{\theta}_i))$. Although the exact distribution of $\bm{y}(k)$ is unknown, it holds that the probability density function of $\bm{y}(k)$ is positive when $\bm{y}(k)$ is in the ellipsoid set $\mathcal{Y}(k|k-1,\bm{\theta}_i)$, while the probability density function is zero when $\bm{y}(k)$ is out of the set. Therefore, we define the likelihood $ L(k,\bm{\theta}_i)$ with the observation $\bm{y}(k)$ as 
\begin{equation}
    \begin{aligned}
        L(k,\bm{\theta}_i) = \left\{ \begin{array}{rl}  \phi_{\text{pdf}}(\bm{y}(k),\bm{\theta}_i), &  \bm{y}(k) \in \mathcal{Y}(k|k-1,\bm{\theta}_i) \\ 0, & \mbox{others} \end{array} \right.
    \end{aligned} 
\end{equation}
where $\phi_{\text{pdf}}(\bm{y}(k),\bm{\theta}_i)>0$. The ergodicity of the process enables that the observation $\bm{y}$ with the true value $\bm{\theta}_i$ converges to $\mathcal{Y}(k|k-1,\bm{\theta}_i)$ and the observation $\bm{y}$ with $\bm{\theta}_j\, (j\neq i)$ diverges from $\mathcal{Y}(k|k-1,\bm{\theta}_i)$. Then, we get the limits
\begin{equation}\label{lim1}
    \begin{aligned}
        \lim_{n \rightarrow \infty} L(k+n-1,\bm{\theta}_i)= \phi_{\text{pdf}}(\bm{y}(k),\bm{\theta}_i),
    \end{aligned} 
\end{equation}
and
\begin{equation}\label{lim2}
    \begin{aligned}
        \lim_{n \rightarrow \infty} L(k+n-1,\bm{\theta}_j)= 0, \quad \forall j \neq i.
    \end{aligned} 
\end{equation}
We define an auxiliary function $\digamma_j(k)$ as
\begin{equation}\label{digamma_1}
\digamma_j(k)= \frac{\pi(k,\bm{\theta}_j)}{\pi(k,\bm{\theta}_i)}, \quad j=1,2, \cdots,s.  
\end{equation}
Substituting the \textit{a-posteriori} probability into the auxiliary function (\ref{digamma_1}), we get
\begin{equation}\label{digamma_2}
\begin{aligned}
\digamma_j(k)&= \frac{L(k,\bm{\theta}_j)\pi(k-1,\bm{\theta}_j)}{L(k,\bm{\theta}_i)\pi(k-1,\bm{\theta}_i)}\\
&=\frac{L(k,\bm{\theta}_j)}{L(k,\bm{\theta}_i)}\digamma_j(k-1).
\end{aligned} 
\end{equation}
According to (\ref{digamma_2}), we have the following equation
\begin{equation}\label{digamma_3}
\frac{\digamma_j(k+n-1)}{\digamma_j(k-1)}=
\frac{\prod_{\tau=k}^{k+n-1}L(\tau,\bm{\theta}_j)}{\prod_{\tau=k}^{k+n-1}L(\tau,\bm{\theta}_i)}. 
\end{equation}
By taking the logarithm of (\ref{digamma_3}), we obtain
\begin{equation}\label{digamma_4}
\begin{aligned}
\ln\left\{ \frac{\digamma_j(k+n-1)}{\digamma_j(k-1)} \right\}&=
\sum_{\tau=k}^{k+n-1}\ln\left\{ L(\tau,\bm{\theta}_j) \right\}\\
&- \sum_{\tau=k}^{k+n-1}\ln\left\{ L(\tau,\bm{\theta}_i) \right\}.
\end{aligned} 
\end{equation}
Taking limits of both sides of (\ref{digamma_4}) and combining \eqref{lim1} and \eqref{lim2}, we obtain
\begin{equation}\label{lim3}
\begin{aligned}
    &\lim_{n \rightarrow \infty }\ln\left\{ \frac{\digamma_j(k+n-1)}{\digamma_j(k-1)} \right\} \\
    &=\lim_{n \rightarrow \infty}\sum_{\tau=k}^{k+n-1}\ln\left\{ L(\tau,\bm{\theta}_j) \right\}
    - \lim_{n \rightarrow \infty}\sum_{\tau=k}^{k+n-1}\ln\left\{ L(\tau,\bm{\theta}_i) \right\} \\
    &=-\infty -\sum_{\tau=k}^{k+n-1} \ln \left\{\phi_{\text{pdf}}(\bm{y}(k),\bm{\theta}_i)\right\} = -\infty.
\end{aligned} 
\end{equation}
Equation \eqref{lim3} can be rewritten as
\begin{equation}\label{lim4}
\lim_{n \rightarrow \infty } \digamma_j(k+n-1) = K e^{-\infty}\digamma_j(k-1)=0, \quad \forall j \neq i ,  
\end{equation}
where $K$ is a constant. Combining \eqref{digamma_1} and \eqref{lim4}, we get
\begin{equation}
\lim_{k \rightarrow \infty }\frac{\pi(k,\bm{\theta}_j)}{\pi(k,\bm{\theta}_i)}=0, \quad \forall j \neq i .  
\end{equation}
Therefore, we have $\pi(k,\bm{\theta}_j) \xrightarrow{} 0$ when $j \neq i$, and $\pi(k,\bm{\theta}_i) \xrightarrow{} 1$ with $\sum_{i=1}^{s}\pi(k,\bm{\theta}_i) = 1$.
\end{proof}

\subsection{Adaptive Robust Control Law}\label{section3_4}
This section derives the adaptive robust control law for the optimization problem in (\ref{prob1}). We derive the control law by aggregating the robust MPC control laws $ \bm{u}(k,\bm{\theta}_i)$ in (\ref{equation_MPC_control_law}) weighted by the posterior probability $\pi(k,\bm{\theta}_i)$ corresponding to the parameter $\bm{\theta}_i$. This method leads to the decoupling between unknown parameter learning and the robust control. \textit{I.e.}, we separate the deriving of control law $ \bm{u}(k,\bm{\theta}_i)$ and the learning of weight $\pi(k,\bm{\theta}_i)$, resulting in computationally simplified implementation of the system control, which is otherwise complicated due to the coupled controlling and parameter learning~\cite{DESHPANDE1973107}. We provide detailed description of the adaptive robust controller as follows.

We define the cost function corresponding to the parameter $\bm{\theta}_i$ at instant $k$ as 
\begin{equation} \label{J_i}
    \begin{aligned}
        J_i(k) =  \left \{ \bm{x}^T(k)\bm{Q}(k)\bm{x}(k)  +\bm{u}^T(k)\bm{R}(k)\bm{u}(k) \right \}_{\bm{\theta}=\bm{\theta}_i}.
    \end{aligned}
\end{equation}
With the \textit{a-posteriori} probability updated by \eqref{PosProb}, we calculate the cost function $J(k)$ as 
\begin{equation}\label{J_k_pi}
    \begin{aligned}
        J(k) &= \bm{x}^T(k)\bm{Q}(k)\bm{x}(k)  +\bm{u}^T(k)\bm{R}(k)\bm{u}(k)\\
        &= \sum_{i=1}^{s}\pi(k,\bm{\theta}_i)J_i(k).
    \end{aligned}
\end{equation}
As we described in Section \ref{section3_2}, we reduce the whole horizon optimization problem in (\ref{prob1}) to an $N_p$-steps predictive optimization problem, where the cost function is $J_k$ in (\ref{J_N}). Using (\ref{J_N}) and (\ref{J_k_pi}), we rewrite the cost function $J_k$ as
\begin{equation}\label{J_N_i}
    \begin{aligned}
        J_k =\sum_{t= k}^{k+N_p-1}  \sum_{i=1}^{s} \pi(t,\bm{\theta}_i)J_i(t).
    \end{aligned}
\end{equation}
Note that the posterior probabilities $\pi(t,\bm{\theta}_i), t=k+1,\cdots,k+N_p-1$ depend on the observations $\bm{y}(t), t=k+1,\cdots,k+N_p-1$, which is not possible to obtain at the instant $k$. In addressing this issue, we use $\pi(k,\bm{\theta}_i)$ to approximate the posterior probabilities $\pi(t,\bm{\theta}_i), t=k+1,\cdots,k+N_p-1$. Therefore, the cost function \eqref{J_N_i} is approximated as
\begin{equation}\label{J_N_i_apr}
    \begin{aligned}
        J_k \approx  \sum_{i=1}^{s}\pi(k,\bm{\theta}_i) \sum_{t= k}^{k+N_p-1}  J_i(t).
    \end{aligned}
\end{equation}
With (\ref{J_N_i_apr}), we obtain the objective of the $N_p$-steps predictive optimization problem as
%\begin{equation}\label{}
%    \begin{aligned}
%        \min_{\{{u}(t)\} \vert_{t=k}^{k+N_p-1}} \max_{ \{\bm{w}(k), \bm{v}(k),\bm{x}(k)\} } J_k
%    \end{aligned}
%\end{equation}
\begin{equation}\label{minmax_J_k}
    \begin{aligned}
        \min_{\{{u}(t)\} \vert_{t=k}^{k+N_p-1}} \max_{ \{\bm{w}(k), \bm{v}(k),\bm{x}(k)\} } \sum_{i=1}^{s}\pi(k,\bm{\theta}_i) \sum_{t= k}^{k+N_p-1}  J_i(t)
    \end{aligned}
\end{equation}
which can be reformulated as a weighted sub-model predictive optimization problem shown below
\begin{equation}\label{minmax_J_i_sum}
    \begin{aligned}
        \sum_{i=1}^{s}\pi(k,\bm{\theta}_i) \left\{\min_{\{{u}(t,\bm{\theta}_i)\} \vert_{t=k}^{k+N_p-1}} \max_{ \{\bm{w}(k), \bm{v}(k),\bm{x}(k)\} } \sum_{t= k}^{k+N_p-1} J_i(t) \right\}.  
    \end{aligned}
\end{equation}
Hence it is reasonable to express the adaptive robust control law as the aggregation of the sub-model control law $\bm{u}(k,\bm{\theta}_i)$ weighted by $\pi(k,\bm{\theta}_i)$, as shown below
\begin{equation}\label{u_adap_robs}
    \bm{u}(k)=\sum_{i=1}^{s}\pi(k,\bm{\theta}_i)\bm{u}(k,\bm{\theta}_i).
\end{equation}

We summarize the main steps in deriving the proposed adaptive robust control law in Table \ref{procedure}. 

\begin{table}[H]
\centering
\caption{\label{procedure}Calculations in the adaptive robust control.}
\resizebox{0.5\textwidth}{!}
{\begin{tabular}{l}
\toprule
\textbf{Initialization:}\\
\tabincell{l}{Initialize the ellipsoid set of system state $E(\bar{\bm{x}}_0,\bm{P}_0)$, the \textit{a-priori}\\ probability $\pi(0,\bm{\theta}_i)$, the parameter vector set $\bm{\Theta}$, the cost function\\ parameter $\bm{Q}(k)$ and $\bm{R}(k)$, the prediction horizon $N_p$, and the \\simulation horizon $N$.}\\
\midrule
\textbf{Computation:}\\
\tabincell{l}{(1) Calculate the robust control law $\bm{u}(k,\bm{\theta}_i), i=1, \cdots, s$ by \\ solving the SDP problem in (\ref{LMI_Theo});}\\
\tabincell{l}{(2) Calculate the \textit{a-posteriori} probability $\pi(k,\theta_i), i=1, \cdots, s$ \\ based on (\ref{PosProb}); }  \\
\tabincell{l}{(3) Calculate the adaptive robust control law $\bm{u}(k)$ based on (\ref{u_adap_robs}) };  \\
\tabincell{l}{(4) Apply the control law $\bm{u}(k)$ to the system and observe the\\ system output $\bm{y}(k+1)$;}\\
\tabincell{l}{(5) Learn the state ellipsoid set $\mathcal{X}(k+1|k+1,\bm{\theta}_i), i=1, \cdots, s$\\ based on the equations (\ref{learning1})-(\ref{learning6});}\\
\tabincell{l}{(6) Repeat steps (1) to (5) for future control iterations.}\\
\bottomrule
\end{tabular}}
\end{table}

\begin{myRemark}
This remark provides a detailed analysis of the computational complexity of our algorithm. Our ellipsoid-set learning algorithm consists of two phases: time prediction and observation update (as described in Lemma~\ref{lemma_sum} and~\ref{lemma_intersect}). The calculations in \eqref{learning1}, \eqref{learning2}, \eqref{learning4}, \eqref{learning5}, \eqref{learning3}, \eqref{learning7}, and \eqref{learning6} in Lemma~\ref{lemma_sum} and~\ref{lemma_intersect} involve multiplication and addition operations of matrices and vectors. Their respective computational complexities are $O(n^2)+O(nr)$, $O(n^3)$, $O(nm)$, $O(n^3)$, $O(n^2m)$, $O(mn)$, and $O(m^3)+O(m^2n)$, calculated according to the Big-Oh notation~\cite{arora2009computational}. The computational complexity in solving the optimal scalar parameter $q(k,\bm{\theta}_i)$ and $q(k,\bm{\theta}_i)$ is $O(n)$ and $O(n^3)+O(n\log(1/\epsilon))$, respectively, where $\epsilon$ is the convergence accuracy and we employ the Newton’s method in seeking the optimal value scalar parameter~\cite{boyd2004convex}. Assuming the state dimension $n$ is much larger than the output number $m$ and input number $r$, the overall computational complexity of the ellipsoid-set learning algorithm is dominated by $O(n^3)$. The complexity of solving the SDP problem \eqref{LMI_Theo} by using interior-point method is $O((N_{\text{dim}})^{3.5})$~\cite{boyd2004convex}, where $N_{\text{dim}}=N_p(n+r)+n+1$ is the dimension of SDP constraint matrix and $N_p$ is the prediction horizon in MPC. 
If we assume the system state dimension $n$ is much larger than the input number $r$, then the dominating computational complexity in solving the SDP comes from $O((n\cdot N_p)^{3.5})$. We perform ellipsoid-set learning and solve the SDP problem $s$ times at each instant. Therefore, the total computational complexity of our adaptive robust control method is $O(s\cdot n^3 + s \cdot (n \cdot N_p)^{3.5})$. We observe that the complexity of our method grows exponentially as the system dimension increases, which may limit the application in large-scale systems. In addressing this issue, sparsity can be adopted in the state matrix, control matrix, and observation matrix to reduce matrix multiplication costs in ellipsoid-set learning. It is also possible to apply the technique proposed by Bertsimas and Brown \cite{bertsimas2007constrained} and develop an inner approximation with second-order cone programming (SOCP) in our robust control method, where solving SOCP is more computationally effective than the SDP for large-scale systems.
\end{myRemark}
%The overall computational complexity of the ellipsoid-set learning algorithm is $O(n^2)+O(nr)+O(n^3)+O(nm)+O(n^3)+O(n^2m)+O(mn)+O(m^3)+O(m^2n)+O(n)+O(n^3)+O(n\log(1/\epsilon))$. Assuming the state dimension $n$ is much larger than the output number $m$ and input number $r$, the dominating computational complexity in ellipsoid-set learning algorithm is attributed by $O(n^3)$. The complexity of solving the SDP problem (54) using the interior-point method is $O((N_{\text{dim}})^{3.5})$~\cite{boyd2004convex}, where $N_{\text{dim}}=N_p(n+r)+n+1$ is the dimension of SDP constraint matrix. Assuming the state dimension $n$ is much larger than the input $r$, the largest computational complexity in solving SDP problem (54) comes from $O((n \cdot N_p)^{3.5})$. 

\begin{myRemark}
The performance of our adaptive robust control will degrade if the true value of uncertain parameters escapes from the finite set of candidate values. In this case, we can use real-time observations to refine the parameter candidate set dynamically. Particularly, we define the ellipsoid set of parameter $\bm{\theta}$ as $E(\hat{\bm{\theta}}(k), \bm{P}_{\bm{\theta}}(k))$, which can be updated by our ellipsoid-set learning algorithm. The center is updated by $\hat{\bm{\theta}}(k+1) = \hat{\bm{\theta}}(k) + \bm{K}_{\bm{\theta}}(k)\{\bm{y}(k+1)-h[\hat{\bm{\theta}}(k), \hat{\bm{x}}(k|k), \bm{u}(k)]\}$, and the matrix is updated by $ \bm{P}_{\bm{\theta}}(k+1) = \beta(k)\{ \left[ \bm{I}-\bm{K}_{\bm{\theta}}(k)\bm{H}(k) \right]            \bm{P}_{\bm{\theta}}(k) \left[ \bm{I}-\bm{K}_{\bm{\theta}}(k)\bm{H}(k) \right]^T + q^{-1}(k)\bm{K}_{\bm{\theta}}(k)\bm{P}_{\delta}(k)     \bm{K}^T_{\bm{\theta}}(k) \}$. The parameter observation equation is $h[\bm{\theta}, \bm{x}(k), \bm{u}(k)]= \bm{C}(k+1,\bm{\theta})[\bm{A}(k,\bm{\theta})\bm{x}(k)+\bm{B}(k,\bm{\theta})\bm{u}(k)]$ and the corresponding observation matrix is $\bm{H}(k) = \frac{\partial h}{\partial\bm{\theta}}\vert_{\hat{\bm{\theta}}(k), \hat{\bm{x}}(k|k), \bm{u}(k)}$. The optimal gain is calculated by $\bm{K}_{\bm{\theta}}(k)=\bm{P}(k)\bm{H}^T(k)\left[ \bm{H}(k) \bm{P}(k)\bm{H}^T(k) +q^{-1}(k)\bm{P}_{\delta}(k) \right] ^{-1}$. The calculation of optimal scalar parameter $q(k)$ and $\beta(k)$ can be referred to Lemma 2. The true parameter is confined in the updated parameter set $E(\hat{\bm{\theta}}(k+1), \bm{P}_{\bm{\theta}}(k+1))$.
We then discretize the set $E(\hat{\bm{\theta}}(k+1), \bm{P}_{\bm{\theta}}(k+1))$ to obtain the renewed parameter space $\bm{\Theta}$.   
\end{myRemark}

\subsection{Extension to Linear Quadratic Tracking Control}\label{section3_5}
The proposed approach can be used to solve linear quadratic regulation problems. In this section, we extend our adaptive robust control approach such that it can also address linear quadratic tracking problems, thus making it a general approach for the control of uncertain linear systems. Below we describe this extension in detail.  

We present the cost function for a linear quadratic tracking problem as
\begin{equation} \label{J_Track_index}
\begin{aligned}
J_T = &\left[\bm{T}\bm{x}(N)-\bm{r}(N)\right]^T\bm{Q}_T(N)\left[\bm{T}\bm{x}(N)-\bm{r}(N)\right]\\
    &+\sum_{k= 0}^{N-1} \left\{ \left[\bm{T}\bm{x}(k)-\bm{r}(k)\right]^T \bm{Q}_T(k)\left[\bm{T}\bm{x}(k)-\bm{r}(k)\right] \right. \\
    &\left. +\bm{u}^T(k)\bm{R}(k)\bm{u}(k) \right\},
\end{aligned}
\end{equation}
where $\bm{r}(k)$ is the reference trajectory, and the matrices $\bm{T}$ and $\bm{Q}_T(k)$ are preset.
We rewrite this cost function with prediction horizon $[k,k+N_p-1]$ as
\begin{equation}\label{J_Track_N}
\begin{aligned}
    &J_k = \sum_{t=k}^{k+N_p-1} \left\{[\bm{T}\bm{x}(t)-\bm{r}(t)]^T\bm{Q}_T(t)[\bm{T}\bm{x}(t)-\bm{r}(t)] \right. \\
        &\quad\quad\quad\quad\quad\quad \left. +\bm{u}^T(t)\bm{R}(t)\bm{u}(t)\right\}.
\end{aligned}
\end{equation}
We rewrite the cost function (\ref{J_Track_N}) as quadratic function of the control sequence 
\begin{equation}\label{J_Track_t_x}
    \begin{aligned}
    J_k = &\bm{U}^T_k\mathscr{B}_k(\bm{\theta}_i)\bm{U}_k +2\bm{b}^T_k(\bm{\theta}_i)\bm{U}_k
     +\bm{W}^T_k\mathscr{C}_k(\bm{\theta}_i)\bm{W}_k\\ &+2\bm{c}^T_k(\bm{\theta}_i) \bm{W}_k +2\bm{U}^T_k\mathscr{D}_k(\bm{\theta}_i)\bm{W}_k +\mathscr{X}_k,
    \end{aligned}
\end{equation}
where the matrices $\mathscr{B}_k(\bm{\theta}_i)$, $\mathscr{C}_k(\bm{\theta}_i)$, and $\mathscr{D}_k(\bm{\theta}_i)$ have been defined in Lemma \ref{lemma3}. We redefine the matrix $\bm{b}_k(\bm{\theta}_i)$, $\bm{c}_k(\bm{\theta}_i)$, and $\mathscr{X}_k$ as 
\begin{equation}
    \bm{b}_k(\bm{\theta}_i) = \sum_{t=k+1}^{k+N_p-1}(\bm{T}\tilde{\bm{B}}^{t-1}_k(\bm{\theta}_i))^T\bm{Q}_T(t)\Delta\bm{x}_r(t),
\end{equation}
\begin{equation}
     \bm{c}_k(\bm{\theta}_i) = \sum_{t=k+1}^{k+N_p-1}(\bm{T}\tilde{\bm{C}}^{t-1}_k(\bm{\theta}_i))^T\bm{Q}_T(t)\Delta\bm{x}_r(t),
\end{equation}
\begin{equation}
\mathscr{X}_k=\sum_{t= k}^{k+N_p-1} \Delta\bm{x}_r^T(t) \bm{Q}_T(t) \Delta\bm{x}_r(t),
\end{equation}
where $\Delta \bm{x}_r(t)$ is defined as
\begin{equation}
\Delta \bm{x}_r(t)= \bm{T}\tilde{\bm{A}}^{t-1}_k(\bm{\theta}_i)\bm{x}(k)-\bm{r}(t),
\end{equation}
and $\bm{Q}(t)$ in matrices $\mathscr{B}_k(\bm{\theta}_i)$, $\mathscr{C}_k(\bm{\theta}_i)$, $\mathscr{D}_k(\bm{\theta}_i)$ is defined as
\begin{equation}
    \bm{Q}(t) = \bm{T}^T\bm{Q}_T(t)\bm{T}.
\end{equation}
Substituting $\bm{x}(k)=\hat{\bm{x}}(k|k,\bm{\theta}_i)+\bm{\eta}(k,\bm{\theta}_i)$ into (\ref{J_Track_t_x}) yields
\begin{equation}\label{J_Track_t_xhat}
    \begin{aligned}
        J_k = &\bm{U}_k^T\mathscr{B}_k(\bm{\theta}_i)\bm{U}_k+2\hat{\bm{b}}_k^T(\bm{\theta}_i)\bm{U}_k +\bm{\xi}_k^T(\bm{\theta}_i)\hat{\mathscr{C}}_k(\bm{\theta}_i)\bm{\xi}_k(\bm{\theta}_i)\\
        &+2\hat{\bm{c}}_k^T(\bm{\theta}_i) \bm{\xi}_k(\bm{\theta}_i)+2\bm{U}_k^T\hat{\mathscr{D}}_k(\bm{\theta}_i)\bm{\xi}_k(\bm{\theta}_i)+\hat{\mathscr{X}}_k,
    \end{aligned}
\end{equation}
where the definitions of matrix $\hat{\mathscr{C}}_k(\bm{\theta}_i)$ and $\hat{\mathscr{D}}_k(\bm{\theta}_i)$ are the same as these presented in Lemma \ref{lemma3}, and the matrix $\hat{\bm{b}}_k(\bm{\theta}_i)$, $\hat{\bm{c}}_k(\bm{\theta}_i)$, and $\hat{\mathscr{X}}_k$ are redefined as 
\begin{equation}
    \hat{\bm{b}}_k(\bm{\theta}_i) = \sum_{t=k+1}^{k+N_p-1}(\bm{T}\tilde{\bm{B}}^{t-1}_k(\bm{\theta}_i))^T\bm{Q}_T(t)\Delta\hat{\bm{x}}_r(t),
\end{equation}
\begin{equation}
\hat{\bm{c}}_k(\bm{\theta}_i)=\left[ \begin{array}{c} \sum_{t=k+1}^{k+N_p-1}(\bm{T}\tilde{\bm{A}}^{t-1}_k(\bm{\theta}_i))^T \bm{Q}_T(t)\Delta\bm{x}_r(t)\\
\sum_{t=k+1}^{k+N_p-1}(\bm{T}\tilde{\bm{C}}^{t-1}_k(\bm{\theta}_i))^T \bm{Q}_T(t)\Delta\hat{\bm{x}}_r(t) \end{array}\right],
\end{equation}
\begin{equation}
    \hat{\mathscr{X}}_k=\sum_{t= k}^{k+N_p-1} \Delta\hat{\bm{x}}_r^T(t) \bm{Q}_T(t) \Delta\hat{\bm{x}}_r(t),
\end{equation}
where $\Delta \hat{\bm{x}}_r(t)$ is
\begin{equation}
\Delta \hat{\bm{x}}_r(t)= \bm{T}\tilde{\bm{A}}^{t-1}_k(\bm{\theta}_i)\hat{\bm{x}}(k)-\bm{r}(t).
\end{equation}
The cost function in (\ref{J_Track_t_xhat}) has the same expression as in (\ref{J_t_xhat}) in Lemma \ref{lemma3}. Therefore, we can apply Theorem \ref{theorem1} to obtain $\bm{u}(k,\bm{\theta}_i)$, and combine it with (\ref{PosProb}) and (\ref{u_adap_robs}) to derive the adaptive robust tracking control law.

\section{Simulations and Results}\label{section5}
In this section, we demonstrate the performance of our adaptive robust control by two numerical simulations. In Section \ref{exp1}, we conduct the linear quadratic regulation for a two-dimensional uncertain system. In Section \ref{exp2}, we consider the linear quadratic tracking control for an uncertain system with three-dimensional states. We also compare our method with the robust control method designed by Bertsimas and Brown~\cite{bertsimas2007constrained}, and the ideal optimal robust control (using ellipsoid-set learning) whose system parameter is known to the controller~\cite{ma2024robust}. 

\subsection{Example 1: adaptive robust linear quadratic regulation control}\label{exp1}
We consider the control problem $(\mathcal{P})$ in (\ref{prob1}) with the following particulars
\begin{equation}
	\begin{aligned}
		\bm{A}(k)=(1+0.2\sin(k)) \left[ \begin{array}{cc} 0.6+\bm{\theta}(1) & 0.7 \\ \bm{\theta}(2) & 0.5+\bm{\theta}(3) \end{array} \right],
	\end{aligned}
\end{equation}
\begin{equation}
	\begin{aligned}
		\bm{B}(k)=\left[ \begin{array}{c} 1+\bm{\theta}(4)  \\ \bm{\theta}(5)  \end{array} \right], 
		\bm{C}(k)=\left[ \begin{array}{cc} 0.2+\bm{\theta}(6) & 1  \end{array} \right],
	\end{aligned}
\end{equation}
\begin{equation}
	\begin{aligned}
		\bm{P}_w(k)=\left[ \begin{array}{cc} (0.1\arctan(k))^2 & 0 \\ 0 & (0.1\arctan(k))^2 \end{array} \right],
	\end{aligned}
\end{equation}
\begin{equation}
	\begin{aligned}
		\bm{P}_v(k)=(0.15\arctan(k))^2 ,
	\end{aligned}
\end{equation}
\begin{equation}
	\begin{aligned}
		\bm{Q}(k)=\left[ \begin{array}{cc} 1 & 0 \\ 0 & 1 \end{array} \right], \quad \bm{R}(k)=1, \quad N = 31.
	\end{aligned}
\end{equation}
The unknown parameter vector $\bm{\theta}=[\bm{\theta}(1),\bm{\theta}(2),\cdots,\bm{\theta}(6)]^T$ takes values from the set $\bm{\Theta}= \{\bm{\theta}_1, \bm{\theta}_2,  \bm{\theta}_3\} $, where $\bm{\theta}_1=[0,0.25,0,0,0.3,0]^T$, $\bm{\theta}_2=[0.4,-0.25,0,-2,0.5,0]^T$, and $\bm{\theta}_3=[0,0.25,0,0.2,0.3,-0.35]^T$. The prior probabilities are set as $\pi(0,\bm{\theta}_1)=\pi(0,\bm{\theta}_2)=\pi(0,\bm{\theta}_3)=1/3$. In this simulation, we assume the true value of $\bm{\theta}$ is $\bm{\theta}_1$, and the initial state is confined in an ellipsoidal set, where the ellipsoid center $\hat{\bm{x}}(0)$ and the shape matrix $\bm{P}(0)$ are set as
\begin{equation}
	\begin{aligned}
		\hat{\bm{x}}(0)= \left[ \begin{array}{c} \hat{x}_1(0)  \\ \hat{x}_2(0)  \end{array} \right]= \left[ \begin{array}{c} 5  \\ -5  \end{array} \right],  \bm{P}(0) = \left[ \begin{array}{cc} 10 & 8 \\ 8 & 10 \end{array} \right].
	\end{aligned}
\end{equation}
We set $N_p$ as $2$. We conduct the simulation in MATLAB, and solve the SDP problem, as described in Theorem~\ref{theorem1}, by utilizing SeDuMi \cite{SeDuMi1} performed in YALMIP.

\begin{figure}[htbp]
	\centering
	\includegraphics[width=0.45\textwidth]{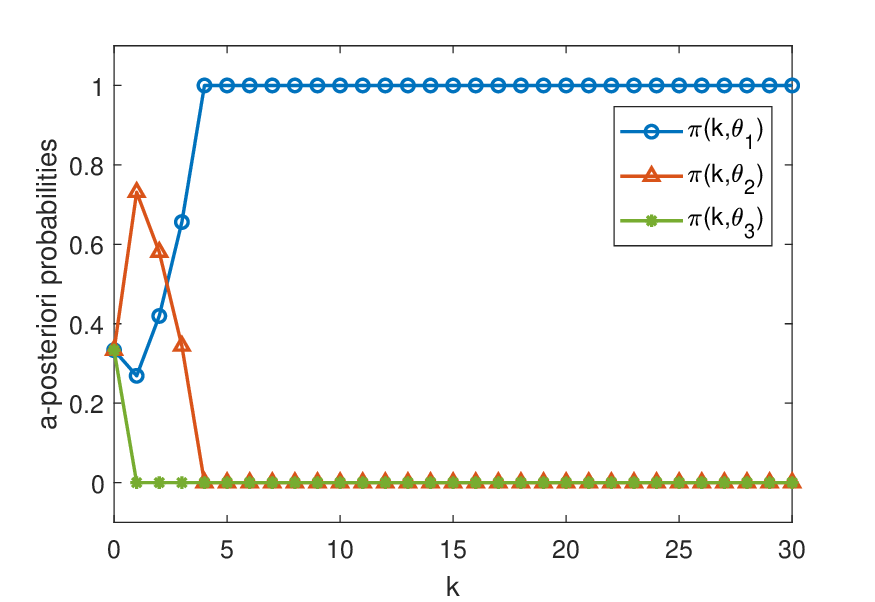}\\
	\caption{The convergence of \textit{a-posteriori} probabilities for candidate system parameters. } \label{prob_converge}
\end{figure}

Fig. \ref{prob_converge} shows the convergence of the weights, or the \textit{a-posteriori} probabilities for the admissible parameter vectors. We observe that the \textit{a-posteriori} probability $\pi(k,\bm{\theta}_1)$ converges to $1$ after $4$ instants, while the \textit{a-posteriori} probabilities of $\pi(k,\bm{\theta}_2)$ and $\pi(k,\bm{\theta}_3)$ converge to $0$. This indicates the system learns and efficiently recognizes the true parameter vector $\bm{\theta}_1$. 

\begin{figure}[htbp]
	\centering
	\includegraphics[width=0.45\textwidth]{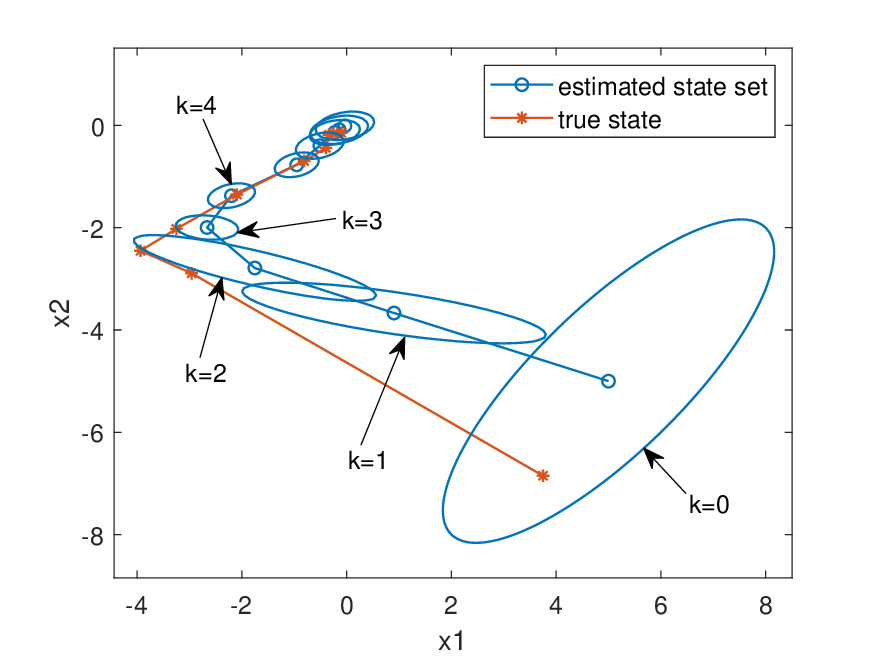}\\
	\caption{The learning progress of ellipsoid set for system states by our adaptive robust control for the system with unknown parameters. } \label{ellip_learn_adaptive}
\end{figure}

\begin{figure}[htbp]
	\centering
	\includegraphics[width=0.45\textwidth]{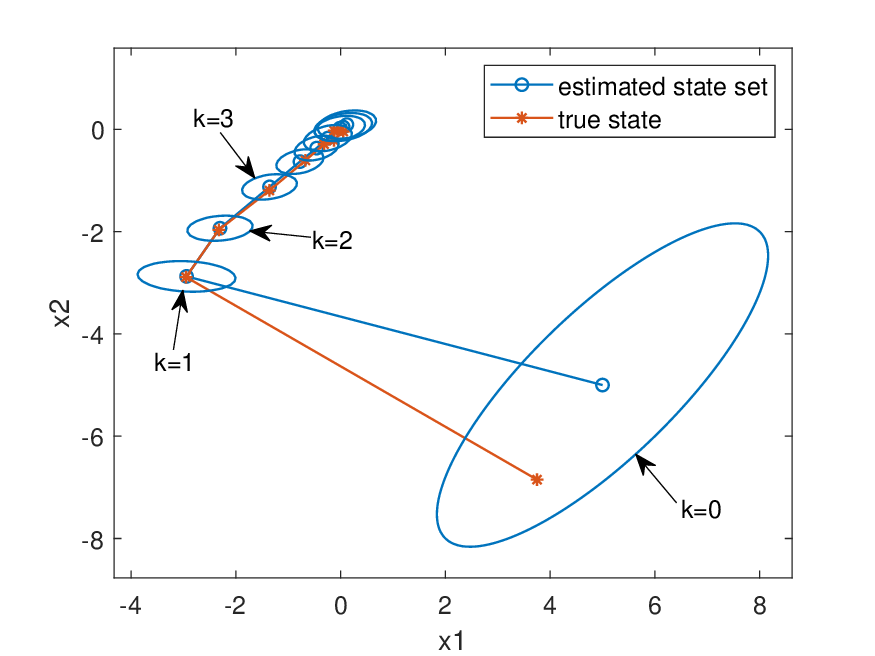}\\
	\caption{The learning progress of ellipsoid set for system states by the ideal optimal robust control for the system with known parameters. } \label{ellip_learn_opt}
\end{figure}

Fig. \ref{ellip_learn_adaptive} and Fig. \ref{ellip_learn_opt} demonstrate the ellipsoid-sets learning progress for the system states by our adaptive robust control and the ideal optimal robust control, respectively. We observe that the ellipsoid sets under both methods shrink over time, and the ellipsoid centers asymptotically converge to the true values, indicating the mitigated uncertainty of the estimated state by the ellipsoid-set learning. In comparison, the ellipsoid set by our adaptive robust control approach converges after $4$ instants, which is slower than the ideal optimal robust control approach (converges within $2$ instants). This is unavoidable because the \textit{a-posteriori} probabilities in our approach need a short time to converge to the right value, since the parameter vector is unknown to the controller. 

\begin{table}[htbp]
\caption{The reduction of state set size utilizing the ellipsoid-set learning.}
\centering
\resizebox{0.5\textwidth}{!}
{\begin{tabular}{ccccccc}
    \toprule
             &  k=0  &  k=1  &  k=2  &  k=3  &  k=4  &  k=5  \\
    \midrule
    $\det(\bm{P}_{\text{OPT}}(k|k))$ & 36 & 0.0850 & 0.0223 & 0.0132 & 0.0102 & 0.0093 \\ 
    $\det(\bm{P}_{\text{ARC}}(k|k))$ & 36 & 0.9292 & 0.8325 & 0.4410 & 0.0130 & 0.0106 \\
    \bottomrule
    \multicolumn{7}{l}{$^*$ $\bm{P}_{\text{OPT}}$ is the shape matrix in the ideal optimal robust control (OPT). }\\
    \multicolumn{7}{l}{$^*$ $\bm{P}_{\text{ARC}}$ is the shape matrix in our adaptive robust control (ARC)}.\\
\end{tabular}}
\label{reductionSet}
\end{table}

In Table \ref{reductionSet}, we use the determinant of the shape matrix to quantify the size of the state ellipsoid set. We observe that the determinant of the shape matrix decreases over time for both OPT and ARC, indicating a progressive reduction in the size of the state set by the ellipsoid-set learning. The set size reduction in ARC is slower than OPT. This is because our ARC needs time to learn the unknown parameter, whereas OPT is the ideal benchmark that assumes full knowledge of the parameter throughout the control.

\begin{figure}[htbp]
\centering
\includegraphics[width=0.45\textwidth]{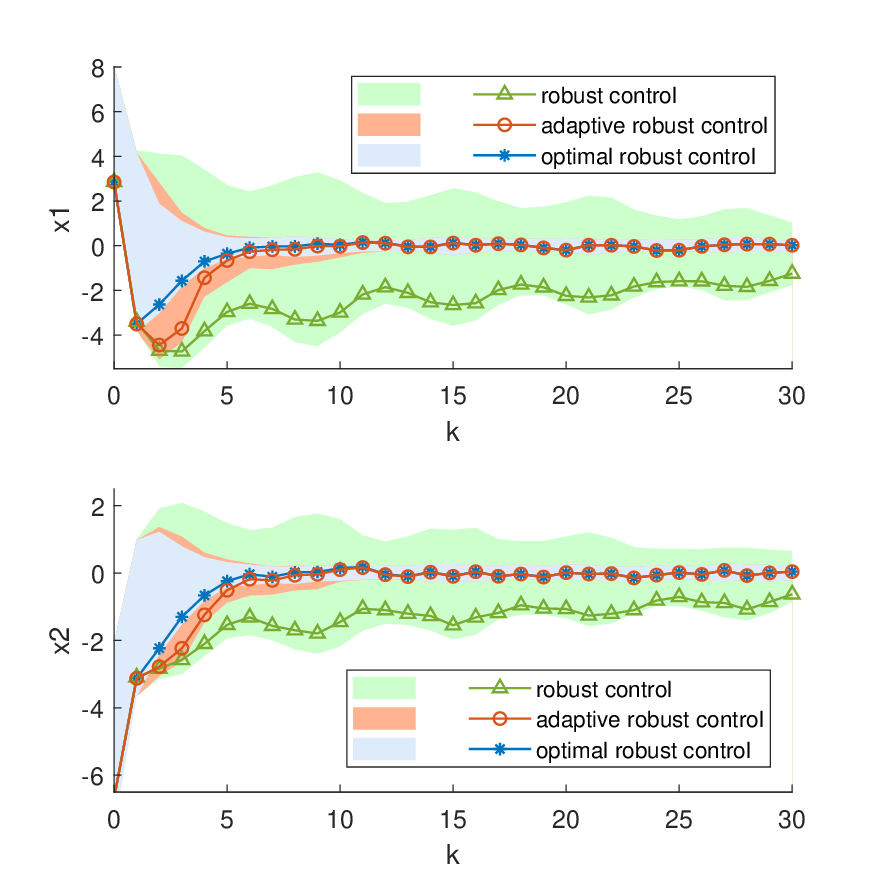}\\
\caption{ Control performances under robust control, our adaptive robust control, and the ideal optimal robust control. The green, red, and blue areas show ranges of the system states during 100 rounds of simulations. The triangle/circle/star-marked lines represent the result from a one-time simulation. } \label{compare1}
\end{figure}

In Fig. \ref{compare1}, we compare the state regulation performances by robust control, our adaptive robust control, and the ideal optimal robust control. We run $100$ rounds of simulations under the three approaches and check whether the system states can be regulated to $0$. We observe that the system states regulation under our approach is better than the robust control, as shown by the closeness of the upper and lower bounds to $0$, implying an enhanced control performance for unknown and bounded noises. Our adaptive robust control follows closely the optimal robust control except for the initial transient period, which reflects the learning toward the unknown system parameters at the beginning of the control.

\begin{table}[]
\caption{Comparison between robust control, adaptive robust control, and optimal robust control.}
\centering
\resizebox{0.5\textwidth}{!}
{\begin{tabular}{ccccccc}
    \toprule
    T         &  5   & 10   & 15   & 20   & 25   & 30 \\ 
    \midrule
    ARC Estimate Error & 1.64 & 1.22 & 1.02 & 0.89 & 0.80 & 0.74 \\ 
    ORC Estimate Error & 0.84 & 0.63 & 0.53 & 0.47 & 0.43 & 0.40 \\
    RC Control Error  & 3.85 & 3.17 & 2.78 & 2.53 & 2.34 & 2.18 \\
    ARC Control Error & 3.59 & 2.67 & 2.21 & 1.93 & 1.74 & 1.59 \\
    ORC Control Error & 3.39 & 2.51 & 2.08 & 1.82 & 0.80 & 1.50 \\
    ORC Estimate Improvement & 94.87\% & 92.83\% & 90.99\% & 88.98\% & 87.38\% & 85.57\%\\ 
    ARC Control Improvement  & 7.12\% & 18.73\% & 25.55\% & 30.70\% & 34.41\% & 36.80\%\\ 
    ORC Control Improvement & 5.84\% & 6.18\% & 6.18\% & 6.17\% & 6.16\% & 6.15\%\\ 
    \bottomrule
    \multicolumn{7}{l}{$^*$ ARC denotes our adaptive robust control, ORC denotes the ideal optimal robust control,}\\
    \multicolumn{7}{l}{ and RC denotes robust control.}
\end{tabular}}
\label{table1}
\end{table}
In Table \ref{table1}, we present the average results of $100$ rounds of simulations. In this table, $T$ is the simulation horizon. ``Estimate Error'' shows the state estimation error varying over time, and is calculated by $\frac{1}{100}\sum_{j=1}^{100} \sum_{k=1}^T \sum_{i=1}^n (\bm{x}_i(k)-\hat{\bm{x}}_i(k))^2_j$. We observe that the estimation errors of ARC and ORC decrease over time with our ellipsoid-set learning algorithm. ``ORC Estimate Improvement'' represents the improvement ratio of ideal benchmark ORC to our ARC. The result shows that the advance of ORC over our ARC is diminishing over time, reflecting that the state estimation performance approaches the one by the ideal benchmark. ``Control Error'' represents the control error over time, calculated by $\frac{1}{100} \sum_{j=1}^{100} \sum_{k=1}^T \sum_{i=1}^n (\bm{x}_i^2(k))_j$. The results show that the control errors under RC, ARC, and ORC decrease over time. ``Control Improvement'' denotes the control performance improvement ratio, where ``ARC Control Improvement''  manifests the improvement of our ARC over RC, and ``ORC Control Improvement'' indicates the advantage of the ideal ORC over our ARC. The result indicated by ``ARC Control Improvement'' shows that ARC brings about $36\%$ improvement compared with RC without ellipsoid-set learning. Note that the improvement increases over time. The result measured by ``ORC Control Improvement'' reflects that the ideal ORC has about $6\%$ improvement compared with our method. Again, we note that this improvement is diminishing over time.

\subsection{Example 2: adaptive robust linear quadratic tracking control}\label{exp2}
We consider the tracking control problem with the following particulars
\begin{equation}
	\begin{aligned}
		\bm{A}(k)=(1+0.5\sin(k)) \left[ \begin{array}{ccc} 0 & 1+\bm{\theta}(1) & 0 \\ 0 & 0 & 1 \\ -0.3 & 0.4 & \bm{\theta}(2) \end{array} \right],
	\end{aligned}
\end{equation}
\begin{equation}
	\begin{aligned}
		\bm{B}(k)=\left[ \begin{array}{c} -0.8\\ 0.7\\ -0.5  \end{array} \right], 
		\bm{C}(k)=\left[ \begin{array}{ccc} 0.5 & 0.8 & 1+\bm{\theta}(3)  \end{array} \right],
	\end{aligned}
\end{equation}
\begin{equation}
    \begin{aligned}
        \bm{P}_w(k)=(0.25\arctan(k))^2 \bm{I}_{3\times 3},
    \end{aligned}
\end{equation}
\begin{equation}
    \bm{P}_v(k)=(0.5\arctan(k))^2,
\end{equation}
\begin{equation}
	\begin{aligned}
		\bm{T}= \left[ \begin{array}{ccc} 1 & 0 & -0.5 \end{array} \right], \bm{Q}_T(k)=10, \bm{R}(k)=1, N = 21.
	\end{aligned}
\end{equation}
The unknown parameter vector $\bm{\theta}=[\bm{\theta}(1),\bm{\theta}(2),\bm{\theta}(3)]^T$ takes values from the set $\bm{\Theta}= \{\bm{\theta}_1, \bm{\theta}_2\} $, where $\bm{\theta}_1=[0,0.2,0]^T$ and $\bm{\theta}_2=[0.2,0.1,-0.5]^T$. The prior probabilities are set as $\pi(0,\bm{\theta}_1)=\pi(0,\bm{\theta}_2)=1/2$. Here we assume the true value of $\bm{\theta}$ is $\bm{\theta}_1$, and the initial state is confined in the ellipsoidal set, where the center $\hat{\bm{x}}(0)$ and shape matrix $\bm{P}(0)$ are set as
\begin{equation}
	\begin{aligned}
		\hat{\bm{x}}(0)= \left[ \begin{array}{c} 10\\ -5\\ 1  \end{array} \right], \quad  \bm{P}(0) = \left[ \begin{array}{ccc} 35 & 20 & 0 \\ 20 & 35 & 0\\ 0 & 0 & 25 \end{array} \right].
	\end{aligned}
\end{equation}
We set $N_p$ as $2$, and generate the reference trajectory by 
\begin{equation}
    \begin{aligned}
        r(k)=\left\{ \begin{array}{rl}
                -4 , & 1\leq k \leq 7 \\ 
                -4+3.5(k-7) , & 7<k\leq9 \\
                3 , & 9<k\leq13 \\
                3-11/6(k-13) , & 13<k\leq16 \\
                -2.5 , & 16<k\leq21
 \end{array}\right. .
    \end{aligned}
\end{equation}

\begin{figure}[htbp]
	\centering
	\includegraphics[width=0.45\textwidth]{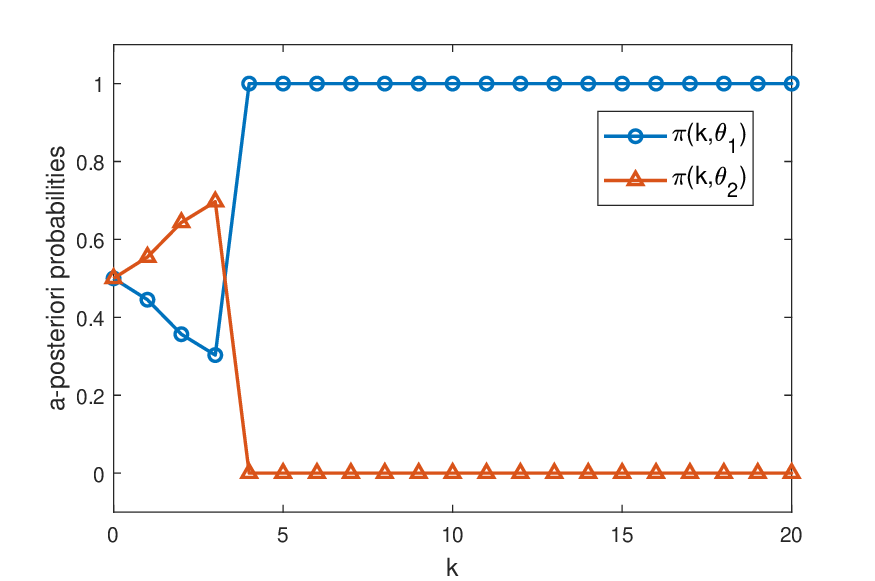}\\
	\caption{The convergence of \textit{a-posteriori} probabilities for candidate parameters under our approach. } \label{prob_converge_2}
\end{figure}

Fig. \ref{prob_converge_2} shows the convergence of the \textit{a-posteriori} probabilities during the parameter estimation of our approach. We observe that the \textit{a-posteriori} probability of the true parameter $\pi(k,\bm{\theta}_1)$ converges to $1$ after $4$ instants and the probability of $\pi(k,\bm{\theta}_2)$ converges to $0$, which implies that the system can correctly learn the parameter vector. 

\begin{figure}[htbp]
	\centering
	\includegraphics[width=0.45\textwidth]{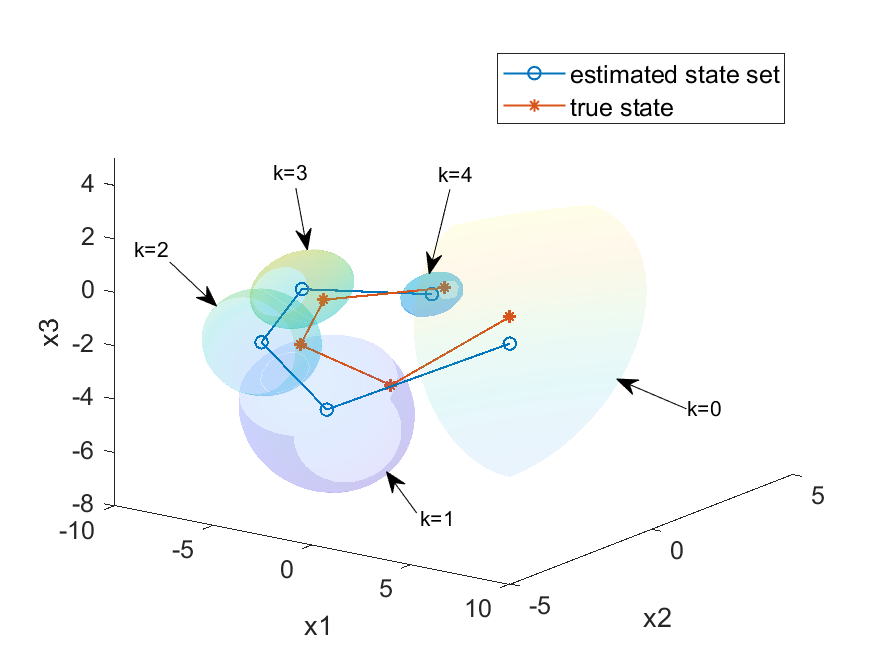}\\
	\caption{The learning progress of ellipsoid set for system states by our adaptive robust control for the system with unknown parameters. } \label{ellip_learn_adaptive_2}
\end{figure}

\begin{figure}[htbp]
	\centering
	\includegraphics[width=0.45\textwidth]{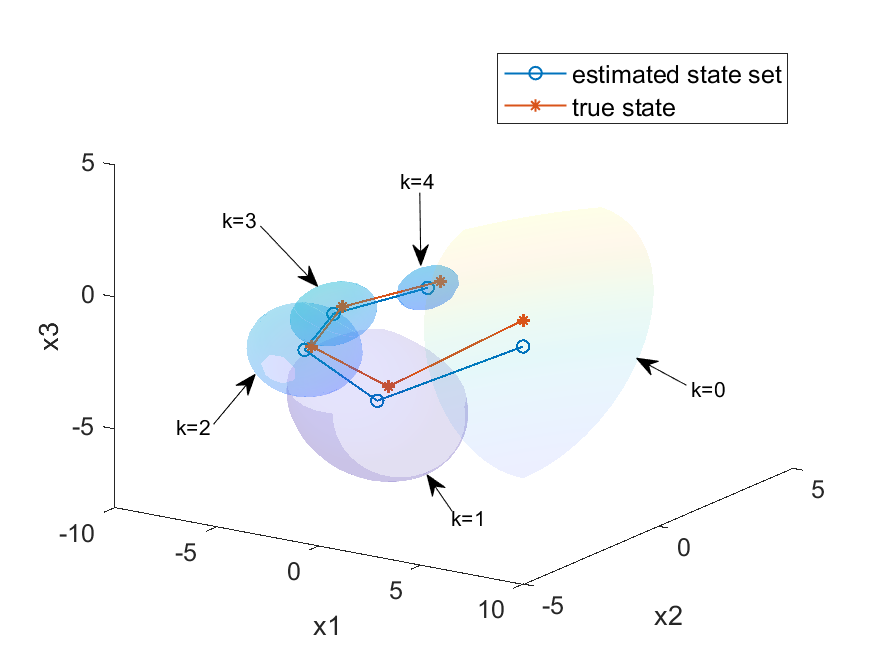}\\
	\caption{The learning progress of ellipsoid set for system states by the ideal optimal robust control for the system with known parameters. } \label{ellip_learn_opt_2}
\end{figure}

Fig. \ref{ellip_learn_adaptive_2} and Fig. \ref{ellip_learn_opt_2} show the ellipsoid-set learning progress for the system states by our approach and the ideal optimal robust control, respectively. We observe that the ellipsoid sets shrink over time, and the ellipsoid centers asymptotically converge to the true values for both approaches, indicating reduced system uncertainty. In comparison, the ellipsoid set by our adaptive robust control approach converges slower than the optimal robust control approach. This is because our approach needs a short time to learn the \textit{a-posteriori} probabilities for the system parameters since they are unknown to the controller. 

\begin{figure}[htbp]
\centering
\includegraphics[width=0.45\textwidth]{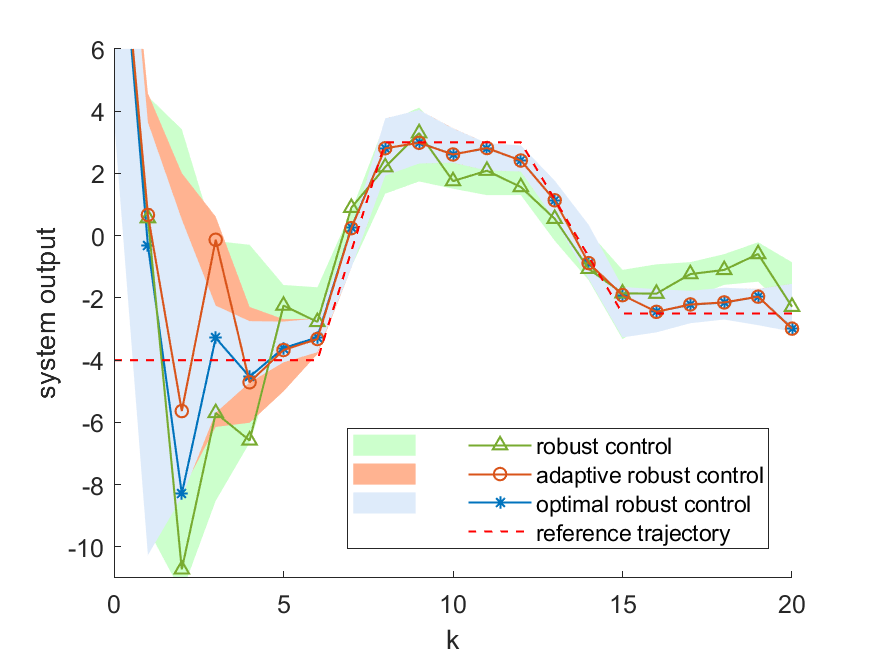}\\
\caption{Control performances under robust control, our adaptive robust control, and the ideal optimal robust control. The green, red, and blue areas show ranges of the system states during 100 rounds of simulations. The triangle/circle/star-marked lines represent the result of a one-time simulation. The red dashed line denotes the reference trajectory.} \label{compare2}
\end{figure}
In Fig. \ref{compare2}, we compare the tracking control performances of robust control, the proposed adaptive robust control, and the ideal optimal robust control. We run $100$ rounds of simulations and check whether the system output can track the reference trajectory. We observe that the system output under our adaptive robust control track the reference trajectory better than the robust control without ellipsoid-set learning, implying an enhanced control performance for the control system with unknown and bounded noises. Our adaptive robust control follows closely the optimal robust control except for the initial transient period, which reflects the essential steps for uncertainty learning toward the unknown system parameters.

\begin{table}[H]
\caption{Comparison between robust control, adaptive robust control, and ideal optimal robust control.}
\centering
\resizebox{0.5\textwidth}{!}
{\begin{tabular}{ccccccc}
    \toprule
    T     &  5   & 10   & 15   & 20    \\ 
    \midrule
    ARC Estimate Error & 2.84 & 2.23 & 1.96 & 1.80 \\ 
    ORC Estimate Error & 2.55 & 2.04 & 1.81 & 1.68 \\
    RC Control Error & 3.42 & 3.36 & 3.24 & 2.92 \\
    ARC Control Error & 3.07 & 2.26 & 1.92 & 1.72 \\
    ORC Control Error & 2.65 & 1.99 & 1.70 & 1.54 \\
    ORC Estimate Improvement & 11.05\% &  9.51\% &  8.23\% &  7.34\% & \\ 
    ARC Control Improvement & 11.24\% & 48.22\% & 69.14\% & 70.34\% & \\ 
    ORC Control Improvement & 15.81\% & 13.95\% & 12.50\% & 11.36\% & \\ 
    \bottomrule
    \multicolumn{5}{l}{$^*$ ARC denotes our adaptive robust control, ORC denotes the ideal optimal}\\
    \multicolumn{5}{l}{ robust control, and RC denotes robust control.}
\end{tabular}}
\label{table2}
\end{table}

In Table \ref{table2}, we present the average results from $100$ rounds of simulations. The state estimation error and estimation improvement ratio are calculated in the same way as in Table~\ref{table1}. The results show that estimation error by the ellipsoid set learning shrinks over time, both for the proposed ARC and the ideal benchmark ORC. In the table, ``ORC Estimate Improvement'' shows that the advantage of the ideal ORC over our ARC is diminishing over time, indicating that the estimated system state of ARC converges to the ones learned by the ideal benchmark. ``Control Error'' represents tracking control error, calculated by $\frac{1}{100} \sum_{j=1}^{100} \sum_{k=1}^T \sum_{i=1}^n ((\bm{T}\bm{x}(k)-\bm{r}(k))^2_i)_j$. The results show that tracking errors under RC, ARC, and ORC decrease over time. ``ARC Control Improvement'' depicts the improvement of our ARC over RC, and ``ORC Control Improvement'' indicates the improvement offered by the ideal ORC over our ARC. We observe that ARC yields about $70\%$ improvement compared to RC without ellipsoid-set learning. Note that such improvement is increasing over time. ``ORC Control Improvement'' indicates that the ideal benchmark ORC has about $11\%$ improvement compared to ARC, while this improvement is diminishing over time.

\section{Conclusion} \label{section6}
This work presents an adaptive robust control for uncertain linear systems whose system parameters and noises are unknown, and the states are not directly measurable. We utilize ellipsoid sets to represent the uncertainties both in the noises and the system states, and develop the ellipsoid-set learning algorithm to estimate the state uncertainty set. This leads to a reduced and concentrated range of worst cases fed to the robust control. Based on the learned ellipsoid set, we derive a robust control law by solving a semidefinite programming (SDP) problem, which reduces the conservativeness typically associated with robust control. To further mitigate parametric uncertainties, we employ the Bayesian posterior probability to learn the unknown parameters. This approach separates the parameter learning from the control law derivation, avoiding coupling between learning and controlling that otherwise complicates the control law derivation. We reflect that solving the SDP problem in the control law derivation can be computationally expensive, particularly for large-scale systems. Therefore, we envision future studies on computational complexity reduction for the proposed robust control, which can limit the overhead for control systems that demand real-time performance.

%We suggest local linearization of the nonlinear systems so that our proposed approach can apply, while the efficiency of the ellipsoid set learning needs to be examined along with the assumptions under the linear conditions.

\section{Appendix} \label{}
\subsection{The ellipsoid set for observation} \label{Appendix_obs}
We rewrite the observation ellipsoid set $\mathcal{O}(k)$ as the same structure as the ellipsoidal set described in the \textit{Notation}. 
With the condition that $m \leq n$ and $\rank(\bm{C}(k,\bm{\theta}_i))=m$, the matrix $\bm{C}(k,\bm{\theta}_i)\in \mathbb{R}^{m\times n}$ has right inverse $\bm{C}^{-1}_r(k,\bm{\theta}_i)$, such that $\bm{C}(k,\bm{\theta}_i)\bm{C}^{-1}_r(k,\bm{\theta}_i)=\bm{I}$. The inverse matrix $\bm{C}^{-1}_r(k,\bm{\theta}_i)$ is given by 
\begin{equation}\label{C_invR}
    \bm{C}^{-1}_r(k,\bm{\theta}_i) = \bm{C}^T(k,\bm{\theta}_i)(\bm{C}(k,\bm{\theta}_i)\bm{C}^T(k,\bm{\theta}_i))^{-1}
\end{equation}
With the inverse matrix (\ref{C_invR}), the observation function described in (\ref{subsys2}) can be rewritten as 
\begin{equation}\label{f_obs1}
    \bm{C}^{-1}_r(k,\bm{\theta}_i)\bm{y}(k)-\bm{x}(k) = \bm{C}^{-1}_r(k,\bm{\theta}_i)\bm{v}(k).
\end{equation}
The observation noise is assumed to be bounded in the ellipsoid set $\bm{v}(k) \in E(\bm{0},\bm{P}_{v}(k))$. Through linear transformation, we have the following ellipsoid set 
\begin{equation}\label{f_obs2}
    \bm{C}^{-1}_r(k,\bm{\theta}_i)\bm{v}(k) \in E(\bm{0},\bm{C}^{-1}_r(k,\bm{\theta}_i)\bm{P}_{v}(k)(\bm{C}^{-1}_r(k,\bm{\theta}_i))^T).
\end{equation}
Substituting Equation (\ref{f_obs1}) into (\ref{f_obs2}), we have 
\begin{equation}\label{f_obs3}
\begin{aligned}
    &(\bm{C}^{-1}_r(k,\bm{\theta}_i)\bm{y}(k)-\bm{x}(k)) \\
    &\in E(\bm{0},\bm{C}^{-1}_r(k,\bm{\theta}_i)\bm{P}_{v}(k)(\bm{C}^{-1}_r(k,\bm{\theta}_i))^T).
\end{aligned}
\end{equation}
To simplify the expression, we make the following definitions
\begin{equation}\label{f_obs4}
    \bm{o}(k)=\bm{C}^{-1}_r(k,\bm{\theta}_i)\bm{y}(k),
\end{equation}
and 
\begin{equation}\label{f_obs5}
    \bm{P_o}(k)=\bm{C}^{-1}_r(k,\bm{\theta}_i)\bm{P}_{v}(k)(\bm{C}^{-1}_r(k,\bm{\theta}_i))^T.
\end{equation}
Combining (\ref{f_obs3}), (\ref{f_obs4}), and (\ref{f_obs5}), we obtain the ellipsoid set for the observation:
\begin{equation}\label{}
    \mathcal{O}(k) = E(\bm{o}(k),\bm{P_o}(k)).
\end{equation}

\subsection{Proof for Lemma 1 and derivation for parameter p} \label{Appendix_Lemma_1}
Consider two ellipsoid sets $E(\bm{a}_1, \bm{Q}_1)$ and $E(\bm{a}_2, \bm{Q}_2)$, where $\bm{Q}_1 $, $\bm{Q}_2 \in \mathbb{R}^{s\times s}$. We define the least size ellipsoid set containing the vector sum of the two ellipsoids as 
\begin{equation}\label{}
    E(\bm{a}, \bm{Q}) \supset E(\bm{a}_1, \bm{Q}_1) \oplus E(\bm{a}_2, \bm{Q}_2).
\end{equation}
According to Chernousko's research \cite{chernousko2000ellipsoidal},  the center $\bm{a}$ and shape matrix $\bm{Q}$ are chosen to be
\begin{equation}\label{E1E2sum}
\begin{aligned}
    & \bm{a} = \bm{a}_1+\bm{a}_2, \\
    & \bm{Q} = (1+p^{-1})\bm{Q}_1+(1+p)\bm{Q}_2 .
\end{aligned}
\end{equation}
Let $E(\bm{a}_1, \bm{Q}_1)=\mathcal{X}'(k|k-1,\bm{\theta}_i)$ and $E(\bm{a}_2, \bm{Q}_2)=\mathcal{W}(k-1)$, we have
$\bm{a}_1=\bm{A}(k-1,\bm{\theta}_i)\hat{\bm{x}}(k-1|k-1,\bm{\theta}_i)+\bm{B}(k-1,\bm{\theta}_i)\bm{u}(k-1)$, $\bm{a}_2=0$,
$\bm{Q}_1=\bm{A}(k-1,\bm{\theta}_i)\bm{P}(k-1|k-1,\bm{\theta}_i)\bm{A}^T(k-1,\bm{\theta}_i)$, and $\bm{Q}_2=\bm{P}_{w}(k-1)$. By substituting $\bm{a}_1$, $\bm{Q}_1$, $\bm{a}_2$, and $\bm{Q}_2$ into (\ref{E1E2sum}), we obtain the ellipsoid $\mathcal{X}(k|k-1,\bm{\theta}_i)$ in Lemma \ref{lemma_sum}.

The parameter $p$ in (\ref{E1E2sum}) is selected by minimizing the trace of $\bm{Q}$. According to the results in \cite{maksarov1996state}, we get the parameter
\begin{equation}\label{p_min}
p= \sqrt{\tr(\bm{Q}_1)/\tr(\bm{Q}_2) }.
\end{equation}
By substituting $\bm{Q}_1=\bm{A}(k-1,\bm{\theta}_i)\bm{P}(k-1|k-1,\bm{\theta}_i)\bm{A}^T(k-1,\bm{\theta}_i)$ and $\bm{Q}_2=\bm{P}_{w}(k-1)$ into (\ref{p_min}), we obtain $p(k,\bm{\theta}_i)$ in (\ref{p_i}).

\subsection{Proof for Lemma 2 and derivation for parameter q} \label{Appendix_Lemma_2}
Consider two ellipsoid sets $E(\bm{a}_1, \bm{Q}_1)$ and $E(\bm{a}_2, \bm{Q}_2)$. We define the least size ellipsoid set containing the intersection of the two ellipsoid sets as
\begin{equation}\label{}
   E(\bm{a}, \bm{Q}) \supset E(\bm{a}_1, \bm{Q}_1) \cap E(\bm{a}_2, \bm{Q}_2).
\end{equation}
According to results in \cite{maksarov1996state}, 
we choose the ellipsoid center $\bm{a}$ and shape matrix $\bm{Q}$ as
\begin{equation}\label{E1E2inter}
\begin{aligned}
    &\bm{a} = \bm{a}_1+\bm{L}(\bm{a}_2-\bm{a}_1), \\
    &\bm{Q} = \beta(q)(\bm{I}-\bm{L})\bm{Q}_1,
\end{aligned}
\end{equation}
where $\bm{L}=\bm{Q}_1(\bm{Q}_1+q^{-1}\bm{Q}_2)^{-1}$, and $\beta(q)=1+q-q(\bm{a}_2-\bm{a}_1)^{T}(q\bm{Q}_1+\bm{Q}_2)^{-1}(\bm{a}_2-\bm{a}_1)$. By substituting $\bm{a}_1$, $\bm{Q}_1$, $\bm{a}_2$, and $\bm{Q}_2$ into (\ref{E1E2inter}), we obtain the set $\mathcal{X}(k|k,\bm{\theta}_i)$ in Lemma \ref{lemma_intersect}.

We select the parameter $q$ in (\ref{E1E2inter}) by minimizing the trace of $\bm{Q}$, which is equivalent to solve the following equation \cite{maksarov1996state}
\begin{equation}\label{q_equation1}
\begin{aligned}
     \frac{\beta'(q)}{\beta(q)} + \sum_{i=1}^{n} \frac{-\frac{d_i(\bm{U}^{-1}\bm{Q}_2\bm{U})}{(\lambda_i(\bm{Q}_1^{-1}\bm{Q}_2)+q)^2}}{\frac{d_i(\bm{U}^{-1}\bm{Q}_2\bm{U})}{\lambda_i(\bm{Q}_1^{-1}\bm{Q}_2)+q}}=0,
\end{aligned}
\end{equation}
where $\beta'(q)=1-(\bm{a}_2-\bm{a}_1)^{T}(q\bm{Q}_1+\bm{Q}_2)^{-1}\bm{Q}_2(q\bm{Q}_1+\bm{Q}_2)^{-1}(\bm{a}_2-\bm{a}_1)$. 
The function \eqref{q_equation1} has no solution when 
\begin{equation}\label{condition1}
(1-(\bm{a}_2-\bm{a}_1)^{T}\bm{Q}_2^{-1}(\bm{a}_2-\bm{a}_1))\tr(\bm{Q}_1)-\tr(\bm{Q}_1\bm{Q}_2^{-1}\bm{Q}_1)>0,
\end{equation}
and $q=0$.
The ellipsoid center and shape matrix of the two ellipsoid sets in Lemma \ref{lemma_intersect} are $\bm{a}_1=\hat{\bm{x}}(k|k-1,\bm{\theta}_i)$, $\bm{Q}_1=\bm{P}(k|k-1,\bm{\theta}_i)$, and $\bm{a}_2=\bm{C}^{-1}_r(k,\bm{\theta}_i)\bm{y}(k)$, $\bm{Q}_2=\bm{C}^{-1}_r(k,\bm{\theta}_i)\bm{P}_{v}(k)(\bm{C}^{-1}_r(k,\bm{\theta}_i))^T$. Substituting them into \eqref{q_equation1} and (\ref{condition1}), we obtain the formulas (\ref{solveq}) and (\ref{cond_nosol}).

\subsection{Proof for Theorem 1} \label{Appendix_Theo2}

According to the definition of $\mathscr{B}_k(\bm{\theta}_i)$ in (\ref{B_math}), $\mathscr{B}_k(\bm{\theta}_i)$ is positive definite since $\bm{Q}(k) \succeq 0$ and $\bm{R}(k) \succ 0$. Therefore, $\mathscr{B}_k^{-1}(\bm{\theta}_i)$ exists and is also symmetric positive definite. Then we can obtain $(\mathscr{B}_k^{-1/2}(\bm{\theta}_i))^T = \mathscr{B}_k^{-1/2}(\bm{\theta}_i)$.
Substituting the control vector in (\ref{U_vecotr}) into the cost function in (\ref{J_t_xhat}) yields
\begin{equation}
    \begin{aligned}
        J_k =  &\bm{z}_k^T\bm{z}_k- \hat{\bm{b}}_k^T(\bm{\theta}_i)\mathscr{B}_k^{-1}(\bm{\theta}_i)\hat{\bm{b}}_k(\bm{\theta}_i)+\bm{\xi}_k^T(\bm{\theta}_i)\hat{\mathscr{C}}_k(\bm{\theta}_i)\bm{\xi}_k(\bm{\theta}_i)\\
        &+2\left(\hat{\bm{c}}_k^T(\bm{\theta}_i) - \hat{\bm{b}}_k(\bm{\theta}_i)^T(\mathscr{B}_k^{-1}(\bm{\theta}_i))^T \hat{\mathscr{D}}_k(\bm{\theta}_i) \right)\bm{\xi}_k(\bm{\theta}_i)\\
        &+2\bm{z}_k^T(\mathscr{B}_k^{-1/2}(\bm{\theta}_i))^T\hat{\mathscr{D}}_k(\bm{\theta}_i)\bm{\xi}_k(\bm{\theta}_i) \\
        &+\hat{\bm{x}}^T(k|k,\bm{\theta}_i)\mathscr{A}_k(\bm{\theta}_i) \hat{\bm{x}}(k|k,\bm{\theta}_i).
    \end{aligned}
\end{equation}
With the definitions in (\ref{h_vector}) and (\ref{F_matrix}), $J_k$ can be rewritten as
\begin{equation}\label{J_const}
    \begin{aligned}
        J_k = \bm{z}_k^T\bm{z}_k + 2\bm{h}_k^T(\bm{\theta}_i)\bm{\xi}_k(\bm{\theta}_i) +2\bm{z}_k^T\bm{F}_k(\bm{\theta}_i)\bm{\xi}_k(\bm{\theta}_i) \\
        +\bm{\xi}_k^T(\bm{\theta}_i) \hat{\mathscr{C}}_k (\bm{\theta}_i)\bm{\xi}_k(\bm{\theta}_i) +const.
    \end{aligned}
\end{equation}
By ignoring the $const$ term in the cost function (\ref{J_const}), we can reform the problem ($\mathcal{P}_i$) in (\ref{prob_pred}) as
\begin{equation}
    \begin{aligned}\label{eq: reformed problem}
        (\mathcal{P}_i): \quad & \min_{\bm{z}_k} \max_{ \{\bm{w}(k), \bm{\eta}(k,\bm{\theta}_i)\} } \left\{ \begin{array}{c}
             \bm{z}_k^T\bm{z}_k + 2\bm{h}_k^T(\bm{\theta}_i)\bm{\xi}_k(\bm{\theta}_i)  \\
             +2\bm{z}_k^T\bm{F}_k(\bm{\theta}_i)\bm{\xi}_k(\bm{\theta}_i)\\ +\bm{\xi}_k^T(\bm{\theta}_i) \hat{\mathscr{C}}(\bm{\theta}_i)\bm{\xi}_k(\bm{\theta}_i)
        \end{array} \right\}\\
          s.t. \quad &\bm{w}^T(k) \bm{P}_w^{-1}(k) \bm{w}(k) \leq 1,  \\
        \quad \quad \quad \quad & \bm{\eta}^T(k,\bm{\theta}_i) \bm{P}^{-1}(k|k,\bm{\theta}_i) \bm{\eta}(k,\bm{\theta}_i) \leq 1.
    \end{aligned}
\end{equation}
We introduce an auxiliary variable $\rho_k$ and rewrite (\ref{eq: reformed problem}) as
\begin{equation}\label{prob_roh}
    \begin{aligned}
        (\mathcal{P}_i): \quad  & \min_{\bm{z}_k} \quad \rho_k \\
        \quad s.t. \quad & \rho_k-\bm{z}_k^T\bm{z}_k - 2\bm{h}_k^T(\bm{\theta}_i)\bm{\xi}_k(\bm{\theta}_i) -2\bm{z}_k^T\bm{F}_k(\bm{\theta}_i)\bm{\xi}_k(\bm{\theta}_i)\\
        \quad \quad \quad \quad &  -\bm{\xi}_k^T(\bm{\theta}_i)\hat{\mathscr{C}}_k(\bm{\theta}_i)\bm{\xi}_k(\bm{\theta}_i) \geq 0 ,\\
        \quad \quad \quad \quad & \bm{W}_k^T \bm{P}_k^W \bm{W}_k \leq N_p , \\
        & \bm{P}_k^W = \diag(\bm{P}_w^{-1}(k), \cdots, \bm{P}_w^{-1}(k+N_p-1)),\\
        \quad \quad \quad \quad & \bm{\eta}^T(k,\bm{\theta}_i) \bm{P}^{-1}(k|k,\bm{\theta}_i) \bm{\eta}(k,\bm{\theta}_i) \leq 1.
    \end{aligned}
\end{equation}
The first constraint in problem (\ref{prob_roh}) can be reformulated as
\begin{equation}\label{constraint1}
    \begin{aligned}
        \left[ \begin{array}{c} 1 \\ \bm{\xi}_k(\bm{\theta}_i) \end{array} \right]^T \bm{M}_{\rho} \left[ \begin{array}{c} 1 \\ \bm{\xi}_k(\bm{\theta}_i) \end{array} \right] \geq 0,
    \end{aligned}
\end{equation}
where
\begin{equation}\label{constr_matrix}
    \begin{aligned}
         \bm{M}_{\rho} = \left[ \begin{array}{cc} \rho_k-\bm{z}_k^T\bm{z}_k & -\bm{h}_k^T(\bm{\theta}_i)-\bm{z}_k^T\bm{F}_k(\bm{\theta}_i) \\ -\bm{h}_k(\bm{\theta}_i)-\bm{F}_k^T(\bm{\theta}_i)\bm{z}_k & -\hat{\mathscr{C}}_k(\bm{\theta}_i) \end{array} \right] .
    \end{aligned}
\end{equation}
The second constraint is rewritten as
\begin{equation}
    \begin{aligned}
        \left[ \begin{array}{c} \bm{\eta}(k,\bm{\theta}_i) \\ \bm{W}_k \end{array} \right]^T \left[ \begin{array}{cc} 0 & 0 \\ 0 & \bm{P}_k^W \end{array} \right] \left[ \begin{array}{c} \bm{\eta}(k,\bm{\theta}_i) \\ \bm{W}_k \end{array} \right] \leq N_p \;.
    \end{aligned}
\end{equation}
With the definition of matrix $\bm{M}_k^{W}$ in (\ref{Matrix_W}), this constraint can be reformed as
\begin{equation}\label{constraint2}
    \begin{aligned}
        \left[ \begin{array}{c} 1 \\ \bm{\xi}_k(\bm{\theta}_i) \end{array} \right]^T \left[ \begin{array}{cc} N_p & 0 \\ 0 & -\bm{M}_k^{W} \end{array} \right] \left[ \begin{array}{c} 1 \\ \bm{\xi}_k(\bm{\theta}_i) \end{array} \right] \geq 0 \;.
    \end{aligned}
\end{equation}
The third constraint can be rewritten as
\begin{equation}
    \begin{aligned}
        \left[ \begin{array}{c} \bm{\eta}(k,\bm{\theta}_i) \\ \bm{W}_k \end{array} \right]^T \left[ \begin{array}{cc} \bm{P}^{-1}(k|k,\bm{\theta}_i) & 0 \\ 0 & 0 \end{array} \right] \left[ \begin{array}{c} \bm{\eta}(k,\bm{\theta}_i) \\ \bm{W}_k \end{array} \right] \leq 1\;.
    \end{aligned}
\end{equation}
With the definition of matrix $\bm{M^{\eta}}(k,\bm{\theta}_i)$ in (\ref{Matrix_eta}), this constraint can be reformed as
\begin{equation}\label{constraint3}
    \begin{aligned}
        \left[ \begin{array}{c} 1 \\ \bm{\xi}_k(\bm{\theta}_i) \end{array} \right]^T \left[ \begin{array}{cc} 1 & 0 \\ 0 & -\bm{M^{\eta}}(k,\bm{\theta}_i) \end{array} \right] \left[ \begin{array}{c} 1 \\ \bm{\xi}_k(\bm{\theta}_i) \end{array} \right] \geq 0
    \end{aligned}
\end{equation}

According to the S-procedure \cite{boyd2004convex}, for all $\bm{\xi}_k(\bm{\theta}_i)$ that satisfies the constraints in (\ref{constraint2}) and (\ref{constraint3}), the constraint in (\ref{constraint1}) also holds if there exist $\tau_1\geq 0$ and $\tau_2\geq 0$ such that
\begin{equation}\label{LMI_spro}
    \begin{aligned}
        &\left[ \begin{array}{cc} \rho_k-\bm{z}_k^T\bm{z}_k & -\bm{h}_k^T(\bm{\theta}_i)-\bm{z}_k^T\bm{F}_k(\bm{\theta}_i) \\ -\bm{h}_k(\bm{\theta}_i)-\bm{F}_k^T(\bm{\theta}_i)\bm{z}_k & -\hat{\mathscr{C}}_k(\bm{\theta}_i) \end{array} \right]\\
        &- \tau_1 \left[ \begin{array}{cc} 1 & 0 \\ 0 & -\bm{M^{\eta}}(k,\bm{\theta}_i) \end{array} \right] - \tau_2 \left[ \begin{array}{cc} N_p & 0 \\ 0 & -\bm{M}_k^{W} \end{array} \right]  \geq 0
    \end{aligned}
\end{equation}
We collect all terms in (\ref{LMI_spro}) and reform them as
\begin{equation}\label{LMI_Schur}
    \begin{aligned}
        &\left[ \begin{array}{cc} \rho_k-\tau_1-\tau_2N_p & -\bm{h}_k^T(\bm{\theta}_i) \\ -\bm{h}_k(\bm{\theta}_i) & \bm{G}_k(\bm{\theta}_i)  \end{array} \right] \\
        &-  \left[ \begin{array}{cc} \bm{z}_k & \bm{F}_k(\bm{\theta}_i)  \end{array} \right]^T \left[ \begin{array}{cc}  \bm{z}_k & \bm{F}_k(\bm{\theta}_i) \end{array} \right]   \geq 0
    \end{aligned}
\end{equation}
According to the Schur complement theorem \cite{boyd2004convex}, (\ref{LMI_Schur}) can be rewritten as (\ref{LMI_Theo}), which completes the proof.

\bibliographystyle{unsrt}
\bibliography{bibfile}

\end{document}